\documentstyle[leqno, amsfonts,amssymb]{article}
\newtheorem{definition}{Definition}[section]
\newtheorem{lemma}[definition]{Lemma}
\newtheorem{theorem}[definition]{Theorem}
\newtheorem{proposition}[definition]{Proposition}
\newtheorem{corollary}[definition]{Corollary}
\newtheorem{remark}[definition]{Remark}

\font\ddpp=msbm10  scaled \magstep 1  
\font\rcmed=msbm9              
\font\rcpeq=msbm6              

\newenvironment{proof}{\noindent{\bf Proof~:}}{\QED\medskip}
\def\QED{\hskip0.1em\hfill\null\ \null\nobreak\hfill
\kern3pt\lower1.8pt\vbox{\hrule\hbox
{\vrule\kern1pt\vbox{\kern1.7pt \hbox{$\scriptstyle
QED$}\kern0.2pt}\kern1pt\vrule}\hrule}}
\def\R{\hbox{\ddpp R}}         
\def\medR{\hbox{\rcmed R}}         
\def\C{{\hbox{\ddpp C}}}       
\def\peqC{{\hbox{\rcpeq C}}}   

\def\frg{{\frak g}}
\def\frh{{\frak h}}
\def\frk{{\frak K}}
\def\frt{{\frak t}}
\def\gc{\frg_\peqC}
\def\Re{{\frak R}{\frak e}\,}
\def\Im{{\frak I}{\frak m}\,}
\def\nilm{\Gamma\backslash G}

\def\db{{\bar{\partial}}}
\def\zzz{{\!\!\!}}

\begin{document}

\title{\bf Hermitian structures on six dimensional nilmanifolds}

\author{\bf Luis Ugarte}

\date{}


\maketitle

\bigskip

{\small {\noindent {\sf Abstract}.- Let $(J,g)$ be a Hermitian
structure on a $6$-dimensional compact nilmanifold $M$ with
invariant complex structure $J$ and compatible metric $g$, which
is not required to be invariant. We show that, up to equivalence
of the complex structure, the strong K\"ahler with torsion
structures $(J,g)$ on $M$ are parametrized by the points in a
subset of the Euclidean space, in particular the region inside a
certain ovaloid corresponds to such structures on the Iwasawa
manifold and the region outside to strong K\"ahler with torsion
structures with nonabelian $J$ on the nilmanifold
$\Gamma\backslash (H^3\times H^3)$, where $H^3$ is the Heisenberg
group. A classification of $6$-dimensional nilmanifolds admitting
balanced Hermitian structures $(J,g)$ is given, and as an
application we classify the nilmanifolds having invariant complex
structures which do not admit Hermitian structure with restricted
holonomy of the Bismut connection contained in SU(3). It is also
shown that on the nilmanifold $\Gamma\backslash (H^3\times H^3)$
the balanced condition is not stable under small deformations.
Finally, we prove that a compact quotient of $H(2,1)\times \medR$,
where $H(2,1)$ is the 5-dimensional generalized Heisenberg group,
is the only 6-dimensional nilmanifold having locally conformal
K\"ahler metrics, and the complex structures underlying such
metrics are all equivalent. Moreover, any invariant locally
conformal K\"ahler metric is a generalized Hopf metric.} }

\medskip

{\small \noindent{\sf Keywords:} {\it Hermitian structure,
K\"ahler with torsion structure, balanced metric, locally
conformal K\"ahler structure, generalized Hopf metric,
nilmanifold}

\medskip

\noindent{\sf MSC 2000:} 53C55; 17B30, 32G05}

\bigskip

\section{Introduction}

Let $(J,g)$ be a Hermitian structure on a manifold $M$, with
fundamental 2-form~$\Omega$ and Lee form $\theta$. The 3-form
$Jd\Omega$ can be identified with the torsion of the Bismut
connection, i.e. the unique Hermitian connection with totally
skew-symmetric torsion~\cite{Bismut,Gau}, and when $Jd\Omega$ is
closed and nonzero (which excludes the K\"ahler case) the
Hermitian structure is called {\it strong K\"ahler with torsion}
(SKT for short)~\cite{AI,FPS}. Such structures arise in a natural
way in physics in the context of supersymmetric $\sigma$-models,
and in general metric connections with totally skew-symmetric
torsion have also applications in type II string theory and black
hole moduli spaces (see~\cite{P} and the references therein).

When the Lee form $\theta$ vanishes identically the Hermitian
structure is called {\it balanced}, and such structures constitute
the class ${\cal W}_3$ in the well-known Gray-Hervella
classification~\cite{GH}. A recent result by Fino and
Grantcharov~\cite{FG} states that for any compact complex manifold
$(M,J)$ with holomorphically trivial canonical bundle, the
existence of a balanced structure $(J,g)$ is a necessary condition
for the existence of a $J$-Hermitian metric on $M$ with vanishing
Ricci tensor of its Bismut connection (see also~\cite{AI,FPS} for
related results).

A Hermitian structure $(J,g)$ is said to be {\it locally conformal
K\"ahler} (LCK for short) if $g$ is conformal to some local
K\"ahler metric in a neighborhood of each point of $M$.  LCK
structures correspond to the Gray-Hervella class ${\cal W}_4$, and
in dimension $\geq 6$ they are characterized by the condition
$d\Omega=\theta\wedge\Omega$.

Let $M$ be a compact Hermitian non-K\"ahler manifold of dimension
$2n \geq 6$. Then the SKT, balanced and LCK conditions are
complementary to each other. In fact, it is well-known that a
K\"ahler metric can be defined as a Hermitian structure in ${\cal
W}_3\cap{\cal W}_4$. Moreover, Alexandrov and Ivanov prove
in~\cite{AI} that $d\Omega\not= \theta\wedge\Omega$ if the
Hermitian structure is SKT (the compactness of $M$ is only needed
here), and a Hermitian structure can only be SKT if $\theta\not=0$
(see also~\cite{FPS}).

In this paper we study SKT, balanced and LCK geometries on
$6$-dimensional compact nilmanifolds $\nilm$ whose underlying
complex structure is invariant, that is, $G$ is a simply-connected
nilpotent Lie group having a discrete subgroup $\Gamma$ such that
the quotient $\nilm$ is compact, and the complex structure on
$\nilm$ stems from a left-invariant one on the Lie group $G$.

We first observe that such study can be reduced to the particular
case when the metric is also invariant. This is shown in~\cite{FG}
for balanced structures using the ``symmetrization'' process,
which is based on a previous idea of Belgun~\cite{Belg}, and we
prove that it also holds for SKT and LCK structures on
nilmanifolds (see Propositions~\ref{symmetrization-SKT}
and~\ref{symmetrization-LCK}). A second reduction comes from the
fact that the study of SKT, balanced and LCK structures can be
carried out up to equivalence of the complex structure. Therefore,
we can restrict our attention to Hermitian structures at the level
of the Lie algebra of $G$ and consider just one representative in
each equivalence class of complex structures. Moreover, in
Section~2 we prove that in dimension six any invariant complex
structure $J$ is equivalent to a complex structure defined by one
of two special types of reduced equations, depending on the
``nilpotency'' of $J$ in the sense of~\cite{CFGU2}.

Salamon proves in~\cite{S} that, up to isomorphism, there are
exactly eighteen 6-dimensional nilpotent Lie algebras admitting
complex structure, which we shall denote here by $\frh_k$ ($1\leq
k\leq 16$), $\frh^-_{19}$ and $\frh^+_{26}$ (see
Theorem~\ref{clasif-complex} for details). For instance, the
nilpotent Lie algebra $\frh_2$ is the Lie algebra of $H^3\times
H^3$, where $H^3$ is the Heisenberg group, $\frh_3$ is the Lie
algebra of $H(2,1)\times \R$, $H(2,1)$ being the 5-dimensional
generalized Heisenberg group, $\frh_5$ is the Lie algebra
underlying the Iwasawa manifold, and $\frh_8$ is the Lie algebra
of $H^3\times \R^3$. In Section~2 it is shown that any complex
structure on $\frh_k$ is nilpotent for $1\leq k\leq 16$, whereas
any complex structure on $\frh^-_{19}$ and $\frh^+_{26}$ is of
nonnilpotent type. Since the structure equations of each one of
these Lie algebras are rational, their corresponding
simply-connected nilpotent Lie groups have a discrete subgroup for
which the quotient is compact~\cite{Mal}.

Fino, Parton and Salamon prove in~\cite{FPS} that a 6-dimensional
compact nilmanifold $\nilm$ admits an invariant SKT structure if
and only if the Lie algebra of $G$ is isomorphic to $\frh_2$,
$\frh_4$, $\frh_5$ or $\frh_8$. In Section~3 we prove that the
same classification is valid if we do not require invariance of
the metric. It is also obtained a more reduced form of the SKT
condition given in~\cite{FPS}, which allows us to show that the
space of SKT structures on a 6-dimensional nilmanifold can be
parametrized, up to equivalence of the complex structure, by the
points in a region of the Euclidean 3-space. More concretely, when
the complex structure is not abelian, there is an ovaloid of
revolution in the Euclidean space such that the region inside
corresponds to SKT structures on the Iwasawa manifold, the region
outside to SKT structures on $\Gamma\backslash (H^3\times H^3)$,
and the points on the ovaloid to SKT structures on the nilmanifold
with underlying Lie algebra $\frh_4$.

A large class of balanced structures is provided by the compact
complex parallelizable manifolds~$M$. In fact, any invariant
compatible metric on~$M$ is balanced~\cite{AG}, and this property
allows us to show in Section~4 that in dimension~$\geq 6$ such
manifolds posses no SKT metrics. We also prove that a compact
nilmanifold $\nilm$ of dimension six admits a balanced metric
compatible with an invariant complex structure if and only if the
Lie algebra of $G$ is isomorphic to $\frh_{19}^-$ or $\frh_{k}$
for some $1\leq k\leq 6$. Fino and Grantcharov construct
in~\cite{FG} a family $J_t$ of invariant complex structures on the
Iwasawa manifold not admitting balanced metrics, except for the
natural complex structure $J_0$. Using their above mentioned
result, this family allows them to conclude that for $t\not=0$ the
complex structure $J_t$ does not admit a Hermitian metric whose
Bismut connection has restricted holonomy in SU(3), providing
counter-examples to a conjecture in~\cite{GIP} as well as the non
stability of this property under small deformations. We show that
the general situation for 6-dimensional compact nilmanifolds
$\nilm$ is the following: there exists an invariant complex
structure on $\nilm$ not admitting a Hermitian metric whose Bismut
connection has restricted holonomy in SU(3) if and only if the Lie
algebra of $G$ is not isomorphic to $\frh_1$, $\frh_6$ or
$\frh_{19}^-$. It is also shown that on the nilmanifold
$\Gamma\backslash (H^3\times H^3)$ the balanced condition is not
stable under small deformations.

Section~5 is devoted to LCK geometry on compact nilmanifolds of
dimension six. We prove that such a nilmanifold $\nilm$ admits an
LCK metric compatible with an invariant complex structure if and
only if the Lie algebra of $G$ is isomorphic to $\frh_1$ or
$\frh_3$, that is, apart from the torus, $\Gamma\backslash
(H(2,1)\times \R)$ is the only 6-dimensional nilmanifold having
LCK structures. It is also shown that the complex structures
underlying such LCK metrics are all equivalent. Moreover, any
invariant LCK metric is a generalized Hopf metric, i.e. the Lee
form is parallel with respect to the Levi-Civita connection. As a
consequence, the only non-toral 5-dimensional nilmanifold
admitting an invariant Sasakian structure is a compact quotient
of~$H(2,1)$.

\section{Invariant complex structures on six dimensional
nilmanifolds}

In this paper we deal with compact complex nilmanifolds
$(M=\nilm,J)$ endowed with an {\it invariant} complex structure
$J$, that is, $G$ is a simply connected nilpotent Lie group and
$\Gamma$ a lattice in $G$ of maximal rank, and $J$ stems from a
left invariant integrable almost complex structure on $G$. Since
the structure is invariant, we can restrict our attention to the
level of the nilpotent Lie algebra $\frg$ of $G$.

Let $\frg$ be a Lie algebra. An endomorphism $J\colon \frg
\longrightarrow \frg$ such that $J^2=-{\rm Id}$ is said to be {\it
integrable} if it satisfies the ``Nijenhuis condition''
$$
[JX,JY]=J[JX,Y]+J[X,JY]+[X,Y],
$$
for any $X,Y\in \frg$. In this case we shall say that $J$ is a
{\it complex structure} on $\frg$.

Let us denote by $\gc$ the complexification of $\frg$ and by
$\gc^*$ its dual, which is canonically identified
to~$(\frg^*)_\peqC$. Given an endomorphism $J\colon \frg
\longrightarrow \frg$ such that $J^2=-{\rm Id}$, there is a
natural bigraduation induced on the complexified exterior algebra
$\bigwedge^* \,\gc^* =\oplus_{p,q} \bigwedge^{p,q}(\frg^*)$, where
the spaces $\bigwedge^{1,0}(\frg^*)$ and
$\bigwedge^{0,1}(\frg^*)$, which we shall also denote by
$\frg^{1,0}$ and $\frg^{0,1}$, are the eigenspaces of the
eigenvalues $\pm i$ of $J$ as an endomorphism of~$\gc^*$,
respectively.

Let $d\colon \bigwedge^* \gc^* \longrightarrow \bigwedge^{*+1}
\gc^*$ be the extension to the complexified exterior algebra of
the usual Chevalley-Eilenberg differential. It is well-known that
$J$ is integrable if and only if $\pi_{0,2} \circ
d\vert_{\frg^{1,0}} \equiv 0$, where $\pi_{p,q}\colon
\bigwedge^{p+q} \gc^* \longrightarrow \bigwedge^{p,q}(\frg^*)$
denotes the canonical projection onto the subspace of forms of
type~$(p,q)$.

\medskip

Next we shall focus on {\it nilpotent} Lie algebras (NLA for
short), that is, the {\it descending central
series}~$\{\frg^k\}_{k\geq 0}$ of $\frg$, which is defined
inductively by
$$
\frg^0 = \frg,\quad\quad \frg^k=[\frg^{k-1},\frg],\quad k\geq 1,
$$
satisfies that $\frg^k=0$ for some $k$. If $s$ is the first
positive integer with this property, then the NLA $\frg$ is said
to be {\it $s$-step} nilpotent.

Salamon proves in~\cite{S} the following equivalent condition for
the integrability of $J$ on a $2n$-dimensional NLA: $J$ is a
complex structure on $\frg$ if and only if $\frg^{1,0}$ has a
basis $\{\omega^j\}_{j=1}^n$ such that $d\omega^1=0$ and
$$
d \omega^{j} \in {\cal I}(\omega^1,\ldots,\omega^{j-1}), \quad
\mbox{ for } j=2,\ldots,n ,
$$
where ${\cal I}(\omega^1,\ldots,\omega^{j-1})$ is the ideal in
$\bigwedge\phantom{\!}^* \,\gc^*$ generated by
$\{\omega^1,\ldots,\omega^{j-1}\}$.

In particular, Salamon's condition in six dimensions is equivalent
to the existence of a basis $\{\omega^j\}_{j=1}^3$ for
$\frg^{1,0}$ satisfying
\begin{equation}\label{ecus}
\cases{\begin{array}{lcl} d\omega^1 \zzz & = &\zzz 0,\\
d\omega^2 \zzz & = &\zzz A_{12}\, \omega^{12} + A_{13}\,
\omega^{13} + A_{1\bar{1}}\, \omega^{1\bar{1}} + A_{1\bar{2}}\,
\omega^{1\bar{2}} + A_{1\bar{3}}\, \omega^{1\bar{3}} \, ,\\
d\omega^3 \zzz & = &\zzz B_{12}\, \omega^{12} + B_{13}\,
\omega^{13} + B_{1\bar{1}}\, \omega^{1\bar{1}} + B_{1\bar{2}}\,
\omega^{1\bar{2}} + B_{1\bar{3}}\, \omega^{1\bar{3}}
\\  \zzz &  &\zzz
+ B_{23}\, \omega^{23} + B_{2\bar{1}}\, \omega^{2\bar{1}} +
B_{2\bar{2}}\, \omega^{2\bar{2}} + B_{2\bar{3}}\,
\omega^{2\bar{3}} \, .
\end{array}}
\end{equation}
for some complex coefficients $A$'s and $B$'s. Here $\omega^{jk}$
(resp. $\omega^{j\overline{k}}$) means the wedge product
$\omega^j\wedge\omega^k$ (resp.
$\omega^j\wedge\omega^{\overline{k}}$), where
$\omega^{\overline{k}}$ indicates the complex conjugation of
$\omega^k$. From now on, we shall use a similar abbreviate
notation for ``basic'' forms of arbitrary bidegree.

\subsection{Reduced form of complex structure equations}

Next we show that there are two special and disjoint types of
complex equations, and that the generic structure
equations~(\ref{ecus}) can always be reduced to one of them,
depending on the ``nilpotency'' of the complex structure.

A complex structure $J$ on a $2n$-dimensional NLA $\frg$ is called
{\it nilpotent} if there is a basis $\{\omega^j\}_{j=1}^n$
for~$\frg^{1,0}$ satisfying $d\omega^1=0$ and
\begin{equation}\label{nilpotent-condition}
d \omega^j \in \bigwedge\phantom{\!\!}^2 \,\langle
\omega^1,\ldots,\omega^{j-1},
\omega^{\overline{1}},\ldots,\omega^{\overline{j-1}} \rangle,
\quad \mbox{ for } j=2,\ldots,n .
\end{equation}
Equivalently~\cite{CFGU2}, the {\it ascending series
$\{\frg^J_l\}_{l\geq 0}$ for $\frg$ adapted to $J$}, which is
defined inductively by $\frg^J_0 =  0$ and
$$\frg^J_l =  \{ X \in {\frg} \, \colon \ [J^k(X), {\frg}]
\subseteq \frg^J_{l-1} \ , k=1,2 \, \} \, ,\quad \mbox{for } l\geq
1,$$ satisfies that $\frg^J_l=\frg$ for some positive integer $l$.

Equations~(\ref{ecus}) encode information not only about the
complex structure $J$, but also about the structure of the
nilpotent Lie algebra $\frg$ itself. Therefore, the coefficients
$A$'s and $B$'s in~(\ref{ecus}) must satisfy those compatibility
conditions imposed by the Jacobi identity of the Lie bracket $[\
,\ ]$ on~$\frg$ (which is equivalent to $d\circ d\equiv 0$), as
well as those conditions ensuring the nilpotency of~$\frg$. For
instance, if $\{Z_j\}_{j=1}^3$ denotes the dual basis of
$\{\omega^j\}_{j=1}^3$, then iterating the bracket $[Z_2,Z_3]$ by
$Z_2$ it is clear that $B_{23}$ must vanish in order to the Lie
algebra $\frg$ be nilpotent. The following result is derived by
imposing these necessary compatibility conditions and it
establishes a first reduction of the generic equations.

\begin{lemma}\label{zerocoeff}\
Let $J$ be a complex structure on an NLA $\frg$ of dimension 6.
\begin{enumerate}
\item[{\rm ({\it a})}] If $J$ is nonnilpotent, then there is a
basis $\{\omega^j\}_{j=1}^3$ for $\frg^{1,0}$
satisfying~$(\ref{ecus})$ with
$A_{1\bar{2}}=B_{1\bar{3}}=B_{23}=B_{2\bar{2}}=B_{2\bar{3}}=0$,
and $A_{1\bar{3}}\not=0$. \item[{\rm ({\it b})}] If $J$ is
nilpotent, then there is a basis $\{\omega^j\}_{j=1}^3$ for
$\frg^{1,0}$ satisfying~$(\ref{ecus})$, where the only
nonvanishing coefficients are among $A_{1\bar{1}}$, $B_{12}$,
$B_{1\bar{1}}$, $B_{1\bar{2}}$, $B_{2\bar{1}}$ and $B_{2\bar{2}}$.
\end{enumerate}
\end{lemma}

The detailed proof of ({\it a}) is given in~\cite[Lemma 2.1,
Proposition 2.2 and Proposition 2.4]{pseudoKahler}. Part ({\it b})
is a direct consequence of~(\ref{nilpotent-condition}).

\medskip

In the following result we give a more reduced form of the
equations for nonnilpotent as well as for nilpotent complex
structures.

\begin{theorem}\label{J-red}
Let $J$ be a complex structure on an NLA $\frg$ of dimension 6.
\begin{enumerate}
\item[{\rm ({\it a})}] If $J$ is nonnilpotent, then there is a
basis $\{\omega^j\}_{j=1}^3$ for $\frg^{1,0}$ such that
\begin{equation}\label{nonnilpotent}
\cases{
\begin{array}{lcl}
d\omega^1 \zzz & = &\zzz 0,\\
d\omega^2 \zzz & = &\zzz E\, \omega^{13} + \omega^{1\bar{3}} \, ,\\
d\omega^3 \zzz & = &\zzz A\, \omega^{1\bar{1}} + ib\,
\omega^{1\bar{2}} - ib\bar{E}\, \omega^{2\bar{1}}  ,
\end{array}
}
\end{equation}
where $A,E\in \C$ with $|E|=1$, and $b\in\R-\{0\}$. \item[{\rm
({\it b})}] If $J$ is nilpotent, then there is a basis
$\{\omega^j\}_{j=1}^3$ for $\frg^{1,0}$ satisfying
\begin{equation}\label{nilpotent}
\cases{
\begin{array}{lcl}
d\omega^1 \zzz & = &\zzz 0,\\
d\omega^2 \zzz & = &\zzz \epsilon\, \omega^{1\bar{1}} \, ,\\
d\omega^3 \zzz & = &\zzz \rho\, \omega^{12} + (1-\epsilon)A\,
\omega^{1\bar{1}} + B\, \omega^{1\bar{2}} + C\, \omega^{2\bar{1}}
+ (1-\epsilon)D\, \omega^{2\bar{2}},
\end{array}
}
\end{equation}
where $A,B,C,D\in \C$, and $\epsilon,\rho \in \{0,1\}$.
\end{enumerate}
\end{theorem}

\begin{proof}
Let us suppose first that $J$ is nonnilpotent. From
Lemma~\ref{zerocoeff}~(a) we have that
$A_{1\bar{2}}=B_{1\bar{3}}=B_{23}=B_{2\bar{2}}=B_{2\bar{3}}=0$ and
$A_{1\bar{3}}\not=0$ in the equations~(\ref{ecus}) for some
(1,0)-basis $\{\omega^j\}$. The remaining coefficients must
guarantee the nilpotency of $\frg$ and the Jacobi identity
$d(d\omega^j)=0$.

Since $0=d(d\omega^2)\wedge\omega^{23\bar{3}} = - A_{1\bar{3}}
\bar{B}_{12}\, \omega^{123\bar{1}\bar{2}\bar{3}}$, the coefficient
$B_{12}$ must be zero. Moreover, from
$0=d(d\omega^2)\wedge\omega^{23\bar{2}} = A_{1\bar{3}}
\bar{B}_{13}\, \omega^{123\bar{1}\bar{2}\bar{3}}$ it follows that
$B_{13}$ also vanishes. Now, the nilpotency of $\frg$ implies that
$A_{12}=0$, because otherwise
$[Z_1,\stackrel{k}{\cdots}[Z_1,[Z_1,Z_2]]\!\cdot\!\cdot\cdot]=(-A_{12})^k
Z_2$ would be a nonzero element in $\frg^k$ for any $k$. In
addition, if we consider the new (1,0)-basis given by
$\tau^1=\omega^1$, $\tau^2=\omega^2$, $\tau^3=\bar{A}_{1\bar{1}}\,
\omega^1 + \bar{A}_{1\bar{3}}\, \omega^3$, then we can suppose
$A_{1\bar{3}}=1$ and $A_{1\bar{1}}=0$.

Therefore, there is a basis $\{\omega^j\}$ of $\frg^{1,0}$
satisfying~$(\ref{ecus})$, where $A_{1\bar{3}}=1$ and the
remaining nonvanishing coefficients are among $A_{13}$,
$B_{1\bar{1}}$, $B_{1\bar{2}}$ and $B_{2\bar{1}}$. Now, since
$$
d(d\omega^2) = (\bar{B}_{1\bar{2}} - A_{13} B_{2\bar{1}})\,
\omega^{12\bar{1}}
$$
and
$$
d(d\omega^3) = (A_{13} B_{2\bar{1}} + B_{1\bar{2}})\,
\omega^{13\bar{1}} - (\bar{A}_{13} B_{1\bar{2}} + B_{2\bar{1}})\,
\omega^{1\bar{1}\bar{3}},
$$
the Jacobi identity implies that the following conditions hold:
$$
A_{13} B_{2\bar{1}} = \bar{B}_{1\bar{2}} = - B_{1\bar{2}}
\quad\mbox{ and }\quad A_{13} \bar{B}_{1\bar{2}} +
\bar{B}_{2\bar{1}}=0.
$$
In particular, $B_{1\bar{2}}=ib$ for some $b\in \R$. Notice that
$b\not=0$ because otherwise $B_{1\bar{2}}$ and $B_{2\bar{1}}$
would be zero and the complex structure $J$ should be nilpotent
(it suffices to interchange $\omega^2$ with $\omega^3$). Finally,
these conditions also imply that $|A_{13}|=1$, so part~({\it a})
of the theorem is proved.

\smallskip

In order to prove part~({\it b}), if $J$ is a nilpotent complex
structure then Lemma~\ref{zerocoeff}~(b) implies the existence of
a (1,0)-basis $\{\omega^j\}$ satisfying~(\ref{ecus}), where all
the coefficients $A$'s vanish except possibly $A_{1\bar{1}}$, and
$B_{13}=B_{1\bar{3}}=B_{23}=B_{2\bar{3}}=0$. Notice that in this
case $(d\circ d) \omega^j=0$, for $j=1,2$. Since $(d\circ d)
\omega^3=B_{2\bar{2}}(-\bar{A}_{1\bar{1}} \omega^{12\bar{1}} +
A_{1\bar{1}} \omega^{1\bar{1}\bar{2}})$, the Jacobi identity of
the Lie bracket implies that $A_{1\bar{1}}B_{2\bar{2}}=0$.

Now, if $A_{1\bar{1}} \not= 0$ then $B_{2\bar{2}}=0$, and we can
suppose $A_{1\bar{1}}=1$ and $B_{1\bar{1}}=0$ after considering
the change of basis $\tau^1=\omega^1$,
$\tau^2=(1/A_{1\bar{1}})\omega^2$ and
$\tau^3=A_{1\bar{1}}\omega^3-B_{1\bar{1}} \omega^2$. Finally,
notice that if the coefficient of $\tau^{12}$ in $d\tau^3$ is
nonzero, then we can normalize it.
\end{proof}

\medskip

For any election of coefficients in the right hand side of
equations~(\ref{nonnilpotent}), resp.~(\ref{nilpotent}), it is
natural to ask whether the resulting equations are ``admissible''
in the sense that there exists a nonnilpotent, resp. nilpotent,
complex structure $J$ on some 6-dimensional NLA $\frg$ having
these equations with respect to some (1,0)-basis. Next we give an
affirmative answer to this question, but first we reformulate it
in more precise terms.

Let $V$ be a real vector space of dimension $2n$, and denote by
$V_\peqC^*$ the dual of the complexification of $V$. Let us fix a
basis $\{\omega^j,\omega^{\overline{j}} \}_{j=1}^n$ for
$V_\peqC^*$, where $\omega^{\overline{j}}$ denotes the complex
conjugate of $\omega^j$. This is equivalent to give an
endomorphism $J\colon V\longrightarrow V$ such that $J^2=-{\rm
Id}_V$, with respect to which the space $V_\peqC^*$ decomposes as
$V_\peqC^* = V^{1,0}\oplus V^{0,1}$, where
$V^{1,0}=\langle\omega^j\rangle$ and $V^{0,1}=\langle
\omega^{\overline{j}} \rangle$ are the eigenspaces of the
eigenvalues $\pm i$ of the extended endomorphism $J\colon
V_\peqC^* \longrightarrow V_\peqC^*$, respectively. Notice that if
$\{X_j,Y_j\}$ is the basis of $V$ dual to the basis
$\{\alpha^j=\frac{1}{2}\Re\omega^j,
\beta^j=\frac{1}{2}\Im\omega^j\}$ of $V^*$, then the endomorphism
$J$ is given by $JX_j=Y_j$, for $j=1,\ldots,n$.

Fixed an $n$-tuple $\mu=(\mu^1,\ldots,\mu^n) \in \bigwedge ^2
V_\peqC^* \times \cdots \times \bigwedge ^2 V_\peqC^*$, we
consider the linear mapping $d_\mu\colon V_\peqC^* \longrightarrow
\bigwedge ^2 V_\peqC^*$ defined by $d_\mu \omega^j=\mu^j$ and
$d_\mu \omega^{\overline{j}} =\overline{\mu^j}$, for
$j=1,\ldots,n$, and we extend it to the complexified exterior
algebra using the formula
$d_\mu(\alpha\wedge\beta)=d_\mu\alpha\wedge\beta + (-1)^{{\rm
deg}\,\alpha} \alpha\wedge d_\mu\beta$, for $\alpha,\beta\in
\bigwedge ^* V_\peqC^*$. Let $[\ ,\ ]_\mu \colon V\times V
\longrightarrow V$ be the bracket on $V$ defined by
$$
[X,Y]_\mu=-\sum_{j=1}^n  \left( \mu^j(X,Y)\, Z_j  +
\overline{\mu^j}(X,Y)\, \bar{Z}_j \right),
$$
for $X,Y\in V$, where $\{Z_j,\bar{Z}_j\}$ is the dual basis of
$\{\omega^j,\omega^{\overline{j}} \}$.

We introduce the following notation: $d_\mu (\mu) \equiv (d_\mu
\mu^1,\ldots,d_\mu \mu^n)$ and $\mu_{0,2} \equiv
(\pi_{0,2}(\mu^1),\ldots,$ $\pi_{0,2}(\mu^n))$, where
$\pi_{0,2}\colon \bigwedge ^2 V_\peqC^* \longrightarrow \bigwedge
^{0,2} (V^*)$ is the canonical projection onto the subspace of
elements of type $(0,2)$.

\begin{lemma}\label{general-converse}
Let $V$ be a real vector space of dimension $2n$, and fix a basis
$\{\omega^j,\omega^{\overline{j}} \}_{j=1}^n$ of $V_\peqC^*$.
Given an $n$-tuple $\mu \in \bigwedge ^2 V_\peqC^* \times \cdots
\times \bigwedge ^2 V_\peqC^*$, we define $J$, $d_\mu$ and $[\ ,\
]_\mu$ as above.
\begin{enumerate}
\item[{\rm ({\it a})}] If $d_\mu(\mu) = 0$, then $\frg_\mu=(V,[\
,\ ]_\mu)$ is a Lie algebra. \item[{\rm ({\it b})}] If in addition
$\mu_{0,2} = 0$, then $J$ is a complex structure on $\frg_\mu$.
\end{enumerate}
\end{lemma}

\begin{proof}
From the definitions we have $\omega^j([X,Y]_\mu)=-d_\mu
\omega^j(X,Y)$. Now ({\it a}) is clear because the bracket~$[\ ,\
]_\mu$ satisfies the Jacobi identity if and only if $d_\mu(d_\mu
\omega^j)=0$ for $j=1,\ldots,n$, that is, $d_\mu(\mu)=0$. To
see~({\it b}), just notice that the Nijenhuis condition is
equivalent to the vanishing of the $(0,2)$-type component in
$d_\mu\omega^j=\mu^j$, for $j=1,\ldots,n$.
\end{proof}

In general, the Lie algebra $\frg_\mu$ may not be nilpotent. For
example, if we consider a 3-tuple
$\mu=(d\omega^1,d\omega^2,d\omega^3)$ given by~(\ref{ecus}) and
satisfying $d_\mu(\mu)=0$, then it determines a Lie algebra
$\frg_\mu$ for which the endomorphism $J$ above is a complex
structure, however $\frg_\mu$ cannot be nilpotent if
$B_{23}\not=0$.

Next we show that for any $\mu$ given by~(\ref{nonnilpotent})
or~(\ref{nilpotent}), we always obtain a nilpotent Lie
algebra~$\frg_\mu$. Thus, the following proposition can be
considered as the converse to Theorem~\ref{J-red}.

\begin{proposition}\label{J-char}
In the conditions of Lemma~$\ref{general-converse}$ we have:
\begin{enumerate}
\item[{\rm ({\it a})}] If $\mu=(0,\ E\, \omega^{13} +
\omega^{1\bar{3}},\ A\, \omega^{1\bar{1}} + ib\, \omega^{1\bar{2}}
- ib\bar{E}\, \omega^{2\bar{1}})$ with $A,E\in \C$, $|E|=1$ and
$b\in \R-\{0\}$, then $\frg_\mu$ is an NLA and $J$ is a
nonnilpotent complex structure on $\frg_\mu$.
 \item[{\rm ({\it b})}] If $\mu=(0,\ \epsilon\, \omega^{1\bar{1}},\
\rho\, \omega^{12} + (1-\epsilon)A\, \omega^{1\bar{1}} + B\,
\omega^{1\bar{2}} + C\, \omega^{2\bar{1}} + (1-\epsilon)D\,
\omega^{2\bar{2}})$ with $A,B,C,D\in \C$ and $\epsilon,\rho\in
\{0,1\}$, then $\frg_\mu$ is an NLA and $J$ is a nilpotent complex
structure on~$\frg_\mu$.
\end{enumerate}
\end{proposition}

\begin{proof}
First, let $\mu$ be given as in ({\it a}). It is easy to check
that $d_\mu(\mu)=0$, so the Jacobi identity holds for the bracket
$[\ ,\ ]_\mu$. In terms of the complex basis $\{Z_j,\bar{Z}_j\}$
dual to $\{\omega^j,\omega^{\overline{j}} \}$, this bracket is
given by
$$
\begin{array}{rcl}
[Z_1,Z_3]_\mu \zzz & = &\zzz - E\, Z_2,\\[1pt]
[Z_1,\bar{Z}_3]_\mu \zzz & = &\zzz - Z_2,\\[1pt]
[Z_1,\bar{Z}_2]_\mu \zzz & = &\zzz -ib\,(Z_3-E\,\bar{Z}_3),\\[1pt]
[Z_1,\bar{Z}_1]_\mu \zzz & = &\zzz - A\, Z_3 + \bar{A}\,\bar{Z}_3,
\end{array}
$$
and their complex conjugates. Therefore, if $E\not=1$ then the
derived algebra $(\frg_\mu)^1=[V,V]_\mu$ is contained in the space
$\langle \Re(Z_2),\Im(Z_2), (1-\bar{E})(Z_3-E\,\bar{Z}_3), i(A\,
Z_3 - \bar{A}\,\bar{Z}_3) \rangle$. Notice that the element
$(1-\bar{E})(Z_3-E\,\bar{Z}_3)$ is in the center of $\frg_\mu$ and
that it is a multiple of $i(A\, Z_3 - \bar{A}\,\bar{Z}_3)$ if and
only if $\bar{A}=AE$. Thus,
$$
\begin{array}{rcl}
(\frg_\mu)^2=[[V,V]_\mu,V]_\mu &\subseteq& \langle \Re(Z_2),\,
\Im(Z_2),\, (1-\bar{E})(Z_3-E\,\bar{Z}_3) \rangle\, ,\\[4pt]
(\frg_\mu)^3=[[[V,V]_\mu,V]_\mu,V]_\mu &\subseteq& \langle
(1-\bar{E})(Z_3-E\,\bar{Z}_3) \rangle\, ,
\end{array}
$$
and $(\frg_\mu)^4=0$, that is, the Lie algebra $\frg_\mu$ is
nilpotent in step $s\leq 4$.

When $E=1$, the elements $i(Z_3 - \bar{Z}_3)$ and $i(A\, Z_3 -
\bar{A}\,\bar{Z}_3)$ of $[V,V]_\mu$ are linearly dependent if and
only if the coefficient $A$ is real. In any case, $i(Z_3 -
\bar{Z}_3)$ is a central element and therefore: if $A\in\R$, then
$\frg_\mu^1=\langle \Re(Z_2),\Im(Z_2), i(Z_3 - \bar{Z}_3)
\rangle$, $\frg_\mu^2=\langle i(Z_3 - \bar{Z}_3) \rangle$ and
$\frg_\mu^3=0$; if $A$ is not real, then $\frg_\mu^1=\langle
\Re(Z_2),\Im(Z_2),\Re(Z_3),\Im(Z_3) \rangle$ and $\frg_\mu$ is
4-step nilpotent.

Finally, the bracket relations above imply that any term in the
ascending series $\{(\frg_\mu)^J_l\}_{l\geq 0}$ adapted to $J$ is
zero, so the complex structure $J$ is nonnilpotent. This completes
the proof of~({\it a})

\smallskip

Now, suppose that $\mu$ is given as in ({\it b}). Since
$d_\mu(\mu)=0$, the bracket $[\ ,\ ]_\mu$ satisfies the Jacobi
identity. The Lie algebra $\frg_\mu=(V,[\ ,\ ]_\mu)$ is nilpotent
in step $s\leq 3$, because $(\frg_\mu)^2=[[V,V]_\mu,V]_\mu
\subseteq \langle \Re Z_3,\Im Z_3 \rangle$, and $\Re Z_3,\Im Z_3$
are central elements of $\frg_\mu$.

The terms in the ascending series $\{(\frg_\mu)^J_l\}_{l\geq 0}$
adapted to $J$ satisfy: $(\frg_\mu)^J_1 \supseteq \langle \Re
Z_3,\Im Z_3=-J(\Re Z_3) \rangle$, $(\frg_\mu)^J_2 \supseteq
\langle \Re Z_2,\Im Z_2=-J(\Re Z_2),\Re Z_3,$ $\Im Z_3=-J(\Re Z_3)
\rangle$, and $(\frg_\mu)^J_3=\frg_\mu$. Therefore, $J$ is a
nilpotent complex structure, and part~({\it b}) of the proposition
is proved.
\end{proof}

\begin{remark}\label{gmu}
{\rm Let us consider a family of $\mu$'s such that $d_\mu(\mu)$
and $\mu_{0,2}$ vanish.
\begin{enumerate}
\item[{\rm ({\it a})}] From Lemma~\ref{general-converse}, we get a
family of Lie algebras $\frg_\mu=(V,[\ ,\ ]_\mu)$ on which the
endomorphism $J\colon V\longrightarrow V$ (which is independent on
$\mu$) is integrable. Let us fix an inner product $\langle\ ,\
\rangle$ on $V$ compatible with $J$ which does not depend on
$\mu$. Now, in the case that $\frg_\mu$ is nilpotent for each
$\mu$, our construction is related to~\cite{L}, where it is
investigated the space of all ``nilpotent'' Lie brackets $[\ ,\
]_\mu$ for which $J$ is integrable and compatible with $\langle\
,\ \rangle$, i.e. $(J,\langle\ ,\ \rangle)$ is a fixed Hermitian
structure on each NLA $\frg_\mu=(V,[\ ,\ ]_\mu)$.
\item[{\rm ({\it b})}] Notice that the Lie algebras $\frg_\mu$
might be nonisomorphic to each other. When $\frg_\mu$ and
$\frg_{\mu'}$ are both isomorphic to a Lie algebra $\frg$, we can
interpret this situation as having two complex structures $J_\mu$
and $J_{\mu'}$ on the same Lie algebra $\frg$.
\end{enumerate}
}
\end{remark}

\subsection{Classification of NLAs admitting
complex structure}

Next we show that a 6-dimensional NLA cannot support nilpotent and
nonnilpotent complex structures at the same time, and then we
classify the NLAs according to the nilpotency of the complex
structures that they admit.

\begin{proposition}\label{1dimensionalcenter}
Let $\frg$ be an NLA of dimension $6$ having a nonnilpotent
complex structure. Then, the center of $\frg$ is $1$-dimensional.
\end{proposition}

\begin{proof}
From Theorem~\ref{J-red}~({\it a}), there is a (1,0)-basis
$\{\omega^j\}_{j=1}^3$ with reduced
equations~(\ref{nonnilpotent}). Then, in terms of its dual basis
$\{Z_j\}$, any central element $T$ of $\frg$ is expressed as
$T=\sum_{j=1}^3 (\lambda_j\, Z_j + \bar{\lambda}_j\, \bar{Z}_j)$,
for some $\lambda_1, \lambda_2,\lambda_3 \in \C$. A direct
calculation shows that $0=[T,Z_3]=-E\lambda_1\, Z_2 -
\bar{\lambda}_1\, \bar{Z}_2$, which implies $\lambda_1=0$.
Moreover,
$$
0=[T,Z_1]=(E \lambda_3 +\bar{\lambda}_3) Z_2 + ib
\bar{\lambda}_2\, Z_3 - ibE\bar{\lambda}_2\, \bar{Z}_3.
$$
Thus $\lambda_2=0$, because $b\not=0$, and $\bar{\lambda}_3=-E
\lambda_3$. Therefore,
$$
T=\lambda_3\, Z_3 - E\lambda_3\, \bar{Z}_3 = \lambda_3\,(Z_3-E\,
\bar{Z}_3).
$$
If $E=1$, then $T=i\lambda(Z_3-\bar{Z}_3)$, $\lambda\in \R$. If
$E\not= 1$, then $T=\lambda(1-\bar{E})(Z_3- E\,\bar{Z}_3)$,
$\lambda\in \R$, because $|E|=1$. Thus,  we conclude that in any
case the center of $\frg$ is 1-dimensional.
\end{proof}

\begin{corollary}\label{or-nilp-or-nonnilp}
Let $\frg$ be a $6$-dimensional NLA admitting complex structures.
Then, all of them are either nilpotent or nonnilpotent.
\end{corollary}

\begin{proof}
If $\frg$ has a nilpotent complex structure $J$, then the first
term $\frg^J_1$ in the ascending series for $\frg$ adapted to $J$
is nonzero. By definition, $\frg^J_1$ is a $J$-invariant ideal of
$\frg$ contained in the center, so if $\frg$ has a nilpotent $J$
then its center is at least 2-dimensional. From
Proposition~\ref{1dimensionalcenter} it follows that $\frg$ has no
nonnilpotent complex structures. Thus, any complex structure on
$\frg$ must be nilpotent.
\end{proof}

\begin{remark}\label{dim10}
{\rm Proposition~\ref{1dimensionalcenter} and
Corollary~\ref{or-nilp-or-nonnilp} do not hold in higher
dimension. In fact, in~\cite{CFGU2} it is given a 10-dimensional
NLA with center of dimension~2 having both nilpotent and
nonnilpotent complex structures. }
\end{remark}

A complex structure $J$ satisfying $[JX,JY]=[X,Y]$, for all
$X,Y\in \frg$, is obviously nilpotent and it is called {\it
abelian}, because $\frg^{1,0}$ is an abelian complex Lie algebra.
It is easily seen that abelian complex structures correspond to
the case $\rho=0$ in the reduced equations~(\ref{nilpotent}).

The following result gives a classification of 6-dimensional NLAs
in terms of the different types of complex structures that they
admit.

\begin{theorem}\label{clasif-complex}
Let $\frg$ be an NLA of dimension~$6$. Then, $\frg$ has a complex
structure if and only if it is isomorphic to one of the following
Lie algebras\footnote{Here we use a mixed notation combining the
structure description of the NLAs as it appears in~\cite{S} and
the notation $\frh_k$ in~\cite{palermo}. For instance,
$\frh_{2}=(0,0,0,0,12,34)$ means that there is a basis
$\{\alpha^j\}_{j=1}^6$ such that the Chevalley-Eilenberg
differential is given by
$d\alpha^1=d\alpha^2=d\alpha^3=d\alpha^4=0$,
$d\alpha^5=\alpha^1\wedge\alpha^2$,
$d\alpha^6=\alpha^3\wedge\alpha^4$; equivalently, the Lie bracket
is given in terms of its dual basis $\{X_j\}_{j=1}^6$ by
$[X_1,X_2]=-X_5$ and $[X_3,X_4]=-X_6$.}
$$
\begin{array}{rcl}
\frh_{1} &\!\!=\!\!& (0,0,0,0,0,0),\\[-2pt]
\frh_{2} &\!\!=\!\!& (0,0,0,0,12,34),\\[-2pt]
\frh_{3} &\!\!=\!\!& (0,0,0,0,0,12+34),\\[-2pt]
\frh_{4} &\!\!=\!\!& (0,0,0,0,12,14+23),\\[-2pt]
\frh_{5} &\!\!=\!\!& (0,0,0,0,13+42,14+23),\\[-2pt]
\frh_{6} &\!\!=\!\!& (0,0,0,0,12,13),\\[-2pt]
\frh_{7} &\!\!=\!\!& (0,0,0,12,13,23),\\[-2pt]
\frh_{8} &\!\!=\!\!& (0,0,0,0,0,12),\\[-2pt]
\frh_{9} &\!\!=\!\!& (0,0,0,0,12,14+25),
\end{array}
\quad\quad\quad
\begin{array}{rcl}
\frh_{10} &\!\!=\!\!& (0,0,0,12,13,14),\\[-2pt]
\frh_{11} &\!\!=\!\!& (0,0,0,12,13,14+23),\\[-2pt]
\frh_{12} &\!\!=\!\!& (0,0,0,12,13,24),\\[-2pt]
\frh_{13} &\!\!=\!\!& (0,0,0,12,13+14,24),\\[-2pt]
\frh_{14} &\!\!=\!\!& (0,0,0,12,14,13+42),\\[-2pt]
\frh_{15} &\!\!=\!\!& (0,0,0,12,13+42,14+23),\\[-2pt]
\frh_{16} &\!\!=\!\!& (0,0,0,12,14,24),\\[-2pt]
\frh^-_{19} &\!\!=\!\!& (0,0,0,12,23,14-35),\\[-2pt]
\frh^+_{26} &\!\!=\!\!& (0,0,12,13,23,14+25).
\end{array}
$$
Moreover:
\begin{enumerate}
\item[{\rm ({\it a})}] Any complex structure on $\frh_{19}^-$ and
$\frh_{26}^+$ is nonnilpotent.
\item[{\rm ({\it b})}] For $1\leq k\leq 16$, any complex structure
on $\frh_{k}$ is nilpotent.
\item[{\rm ({\it c})}] Any complex structure on $\frh_1$,
$\frh_{3}$, $\frh_{8}$ and $\frh_{9}$ is abelian.
\item[{\rm ({\it d})}] There exist both abelian and nonabelian
nilpotent complex structures on $\frh_{2}$, $\frh_{4}$, $\frh_{5}$
and $\frh_{15}$.
\item[{\rm ({\it e})}] Any complex structure on $\frh_{6}$,
$\frh_{7}$, $\frh_{10}$, $\frh_{11}$, $\frh_{12}$, $\frh_{13}$,
$\frh_{14}$ and $\frh_{16}$ is not abelian.
\end{enumerate}
\end{theorem}

\begin{proof}
Salamon proves in~\cite{S} that $\frg$ has a complex structure $J$
if and only if it is isomorphic to one of the Lie algebras
appearing in the list above. Now, using
Proposition~\ref{1dimensionalcenter} we have that a
nonnilpotent~$J$ can only live on $\frh_{19}^-$ or $\frh_{26}^+$,
because the center of these NLAs is 1-dimensional.
Corollary~\ref{or-nilp-or-nonnilp} implies that any $J$ on
$\frh_{19}^-$ and $\frh_{26}^+$ is nonnilpotent, and ({\it a}) is
proved.

In \cite{palermo} it is shown that if $J$ is nilpotent then $\frg$
must be isomorphic to $\frh_{k}$ for some $1\leq k\leq 16$. By
Corollary~\ref{or-nilp-or-nonnilp}, any complex structure on
$\frh_k$, $1\leq k\leq 16$, is nilpotent, so ({\it b}) is proved.

Since $\frh_{3}$ and $\frh_{8}$ have first Betti number $\dim
(\frh/[\frh,\frh])$ equal to $5$, any complex structure must be
abelian. On the other hand, since the Lie algebra $\frh_{9}$ is
3-step and its complex structures are all nilpotent, the
coefficient $\epsilon$ in~(\ref{nilpotent}) must be equal to~1.
Therefore, $\rho=0$ because the first Betti number of $\frh_{9}$
is equal to~4, so ({\it c}) is proved.

In~\cite{abelian} it is proved that a 6-dimensional nilpotent Lie
algebra admits an abelian $J$ if and only if it is isomorphic to
$\frh_k$, for $k=1,2,3,4,5,8,9$ or $15$. This proves ({\it e}).

Finally, to see ({\it d}) we observe that the equations
$$
d\omega^1=d\omega^2=0,\quad d\omega^3=\omega^{12} +
\omega^{1\bar{2}} +C\,\omega^{2\bar{1}},
$$
define, in the sense of Proposition~\ref{J-char}~({\it b}), a
nilpotent complex structure on $\frh_2$ for $C=1$, and on $\frh_4$
for $C=2$. On the other hand, the equations
$$
d\omega^1=0,\quad d\omega^2= \epsilon\, \omega^{1\bar{1}} ,\quad
d\omega^3=\omega^{12},
$$
define a nilpotent complex structure on $\frh_5$ for $\epsilon=0$,
and on $\frh_{15}$ for $\epsilon=1$. Since in each case the
coefficient of $\omega^{12}$ in $d\omega^3$ is nonzero, the
complex structures are not abelian. This, together with the fact
that $\frh_k$ has abelian complex structures for $k=2,4,5$ and
$15$, proves ({\it d}) and so the proof of the theorem is
completed.
\end{proof}

\begin{remark}\label{complex-par}
{\rm If $\frg$ is a {\it complex} Lie algebra, then its canonical
complex structure $J$ satisfies $[JX,Y]=J[X,Y]$, for all $X,Y\in
\frg$. Any complex structure $J$ on an NLA $\frg$ satisfying this
condition is obviously nilpotent. Moreover $d(\frg^{1,0})\subset
\bigwedge^{2,0}(\frg^*)$, so in dimension~6 the corresponding
equations are of the form~(\ref{nilpotent}) with $\rho=0,1$ and
all the remaining coefficients equal to zero. Therefore, these
complex structures only live on the abelian Lie algebra $\frh_1$
and on the Lie algebra $\frh_5$ underlying the Iwasawa manifold.
We shall refer to them as {\it complex parallelizable} structures,
because the corresponding complex nilmanifolds posses three
holomorphic 1-forms which are linearly independent at each point.
}
\end{remark}

\begin{remark}\label{deformations-2-step}
{\rm The deformation of abelian invariant complex structures on
2-step nilmanifolds is studied in~\cite{MPPS}, where it is proved
that the Kuranishi process preserves the invariance of the
deformed complex structures, at least for small deformations.
Conditions under which the deformed structures remain abelian are
also investigated there. In this context, it follows from
Theorem~\ref{clasif-complex} that in dimension~6 all the complex
structures obtained by such small deformations are always of
nilpotent type.}
\end{remark}

As a consequence of Theorem~\ref{clasif-complex} we find reduced
complex structure equations for the Lie algebras $\frh_{19}^-$ and
$\frh_{26}^+$.

\begin{proposition}\label{h19-h26}
For any complex structure on $\frh_{19}^-$ (resp. on
$\frh_{26}^+$), there is a $(1,0)$-basis
satisfying~$(\ref{nonnilpotent})$ with $\bar{A}=AE$ (resp.
$\bar{A}\not= AE$).
\end{proposition}

\begin{proof}
Since any complex structure on $\frh_{19}^-$ and $\frh_{26}^+$ is
nonnilpotent, there exist a $(1,0)$-basis
satisfying~$(\ref{nonnilpotent})$. But an NLA $\frg$ defined
by~(\ref{nonnilpotent}) is isomorphic to $\frh_{19}^-$ if and only
if its first Betti number $\dim (\frg/[\frg,\frg])$ is equal to
$3$, which is equivalent to the closedness of the real 1-form
$i(1-\bar{E})(E\omega^3+\omega^{\bar{3}})$. This latter condition
is satisfied if and only if $\bar{A}=AE$.
\end{proof}

\subsection{Complex structure equations on $2$-step NLAs}

Here we shall arrive at more reduced equations which describe any
complex structure on each $2$-step~NLA.

\begin{proposition}\label{epsilon}
Let $\frg$ be a $6$-dimensional NLA endowed with a nilpotent
complex structure~$J$. Then, the coefficient $\epsilon$ vanishes
in the reduced equations~$(\ref{nilpotent})$ corresponding to $J$
if and only if the Lie algebra $\frg$ is nilpotent in step $s\leq
2$ and its first Betti number is $\geq 4$. In this case, $\frg$
must be isomorphic to $\frh_8$ or $\frh_k$ for some $1\leq k\leq
6$.
\end{proposition}

\begin{proof}
It is clear that $\epsilon=0$ in~(\ref{nilpotent}) implies that
$\frg$ is nilpotent in step $\leq 2$ and $\dim (\frg/[\frg,\frg])
\geq 4$.

Suppose that the Lie algebra $\frg$ has first Betti number $\geq
4$ and it is nilpotent in step $s\leq 2$. Let~(\ref{nilpotent}) be
equations corresponding to $J$ on $\frg$, and suppose that
$\epsilon=1$. First, the coefficient $\rho$ must vanish, because
otherwise the first Betti number would be 3. Moreover, if $B$ and
$C$ are not both zero, then $BC\not=0$ in order to have first
Betti number at least 4. Now, if $\{Z_j\}$ is the dual basis of
$\{\omega^j\}$, then the element $\Im Z_2 \in [\frg,\frg]$
satisfies $[\Im Z_2,\frg]\not=0$, that is, the Lie algebra is not
nilpotent in step $s\leq 2$. Therefore, if $\epsilon=1$ then
$B=C=0$, but in such case we can choose $\epsilon=0$ after
interchanging $\omega^2$ with $\omega^3$.

Finally, if $\frg$ has first Betti number $\geq 4$, then
Theorem~\ref{clasif-complex} implies that $\frg$ cannot be
isomorphic to $\frh_7$ or $\frh_k$ for any $k\geq 10$. On the
other hand, $\frh_9$ is 3-step nilpotent, so $\frg$ cannot be
isomorphic to $\frh_9$ if $\epsilon=0$.
\end{proof}

The following lemma provides a further reduction of the equations
on $2$-step NLAs.

\begin{lemma}\label{epsilonzero-new}
Let $J$ be a complex structure on a $2$-step NLA $\frg$ of
dimension~$6$ with first Betti number~$\geq 4$. If $J$ is not
complex parallelizable, then there is a basis
$\{\omega^j\}_{j=1}^3$ of $\frg^{1,0}$ such that
\begin{equation}\label{epsilonzero-red-new}
\cases{d \omega^1=d\omega^2=0,\cr d\omega^3=\rho\, \omega^{12} +
\omega^{1\bar{1}} + B\,\omega^{1\bar{2}} +
D\,\omega^{2\bar{2}},\cr}
\end{equation}
where $B,D\in \C$, and $\rho= 0,1$.
\end{lemma}

\begin{proof}
First, by the preceding proposition we can suppose $\epsilon=0$ in
the reduced equations~$(\ref{nilpotent})$ corresponding to $J$.
Next, we distinguish several cases depending on the vanishing of
the coefficients $A$ and~$D$.

If $A\not=0$, then we consider the change of basis given by
$\omega^1=\omega'^1-C\,\omega'^2$, $\omega^2=A\,\omega'^2$,
$\omega^3=A\, \omega'^3$. It is easy to check that with respect to
the new $(1,0)$-basis $\{\omega'^j\}$ the equations become
\begin{equation}\label{epsilonzero-red-1}
d \omega'^1=d\omega'^2=0,\quad\quad d\omega'^3=\rho\, \omega'^{12}
+ \omega'^{1\bar{1}} + B'\,\omega'^{1\bar{2}} +
D'\,\omega'^{2\bar{2}},
\end{equation}
where $B'=(\bar{A}B-A\bar{C})/A$, and $D'=\bar{A}(AD-BC)/A$.

The case $D\not=0$ is reduced to the previous one if we
interchange $\omega^1$ with $\omega^2$, and change the sign
of~$\omega^3$. Notice that in this case we
get~$(\ref{epsilonzero-red-1})$ with $B'=(\bar{B}D-C\bar{D})/D$,
and $D'=\bar{D}(AD-BC)/D$.

Let us suppose $A=D=0$ in equations~(\ref{nilpotent}). The change
of basis given by $\omega^1 = \omega''^1 + \omega''^2$, $\omega^2
= \omega''^1 - \omega''^2$, $\omega^3 = -2 \omega''^3$,
transforms~(\ref{nilpotent}) into equations of the form
$$
d \omega''^1=d\omega''^2=0,\quad\quad d\omega''^3=\rho\,
\omega''^{12} + A''\, \omega''^{1\bar{1}} +
B''\,\omega''^{1\bar{2}} + C''\,\omega''^{2\bar{1}} +
D''\,\omega''^{2\bar{2}},
$$
where $D''=-A''=(B+C)/2$, and $B''=-C''=(B-C)/2$. Therefore, if
$B+C\not=0$, then we can again reduce these equations to the
form~$(\ref{epsilonzero-red-1})$ with $B'=(|B|^2-|C|^2)/(B+C)$ and
$D'=-BC(\bar{B}+\bar{C})/(B+C)$. Finally, if $A=D=B+C=0$ then
using the change of basis given by $\omega'^1=\omega''^1 +
i\,\omega''^2$, $\omega'^2= i\,\omega''^1 + \omega''^2$ and
$\omega'^3= 2\,\omega''^3$, we arrive at equations of the form
$$d \omega'^1=d\omega'^2=0,\quad\quad
d\omega'^3=\rho\, \omega'^{12} + iB (\omega'^{1\bar{1}} -
\omega'^{2\bar{2}}).
$$
Now, if $J$ is not complex parallelizable then the coefficient
$B\not=0$ and we can apply the argument used in the case
``$A\not=0$'' above to get equations of the
form~(\ref{epsilonzero-red-1}), with $B'=0$ and $D'=-|B|^2$.
\end{proof}

\begin{lemma}\label{center}
Let $J$ be a complex structure on an NLA $\frg$ with reduced
equations~$(\ref{epsilonzero-red-new})$. Then, the dimension of
the center of $\frg$ is $\geq 3$ if and only if $|B|=\rho$ and
$D=0$.
\end{lemma}

\begin{proof}
Let $\{Z_j\}$ be the dual basis of $\{\omega^j\}$. From
equations~(\ref{epsilonzero-red-new}) it is clear that $\Re Z_3$
and $\Im Z_3$ belong to the center of $\frg$. Now, if
$T=\sum_{j=1}^2 (\lambda_j\, Z_j + \bar{\lambda}_j\, \bar{Z}_j)$
is a central element in $\frg$ for some $(\lambda_1, \lambda_2)
\in \C^2$, then the condition
$$
0=[T,Z_1]=(\rho\lambda_2 + \bar{\lambda}_1 + B\, \bar{\lambda}_2)
Z_3 - \bar{\lambda}_1\, \bar{Z}_3
$$
implies that $\lambda_1$ must be zero. In addition, there is a
solution $\lambda_2\not=0$ of the equation $B\,\bar{\lambda}_2 +
\rho\lambda_2=0$ if and only if $|B|=\rho$. Moreover, the
condition $0=[T,Z_2]= D\, \bar{\lambda}_2\, Z_3 - \bar{D}\,
\bar{\lambda}_2\, \bar{Z}_3$ implies that $D=0$ if $\lambda_2\not=
0$. Therefore, there is an element $T$ in the center of $\frg$
such that $\{\Re Z_3,\Im Z_3,T\}$ are linearly independent if and
only if $|B|=\rho$ and $D=0$.
\end{proof}

We finish this section with a general result showing which are, in
the sense of Proposition~\ref{J-char}~({\it b}), the NLAs
underlying the reduced equations~(\ref{epsilonzero-red-new}) in
terms of the coefficients $\rho$, $B$ and $D$.

\begin{proposition}\label{general}
Let $J$ be a complex structure on an NLA $\frg$ given
by~$(\ref{epsilonzero-red-new})$, and let us denote $x=\Re D$ and
$y=\Im D$. Then:
\begin{enumerate}
\item[{\rm (i)}] If $|B|=\rho$, then the Lie algebra $\frg$ is
isomorphic to
\begin{enumerate}
\item[{\rm (i.1)}] $\frh_2$, for $y\not=0$; \item[{\rm (i.2)}]
$\frh_3$, for $\rho=y=0$ and $x\not=0$; \item[{\rm (i.3)}]
$\frh_4$, for $\rho=1$, $y=0$ and $x\not=0$; \item[{\rm (i.4)}]
$\frh_6$, for $\rho=1$ and $x=y=0$; \item[{\rm (i.5)}] $\frh_8$,
for $\rho=x=y=0$.
\end{enumerate}
\item[{\rm (ii)}] If $|B|\not=\rho$, then the Lie algebra $\frg$
is isomorphic to
\begin{enumerate}
\item[{\rm (ii.1)}] $\frh_2$, for $4y^2 >
(\rho-|B|^2)(4x+\rho-|B|^2)$; \item[{\rm (ii.2)}] $\frh_4$, for
$4y^2 = (\rho-|B|^2)(4x+\rho-|B|^2)$; \item[{\rm (ii.3)}]
$\frh_5$, for $4y^2 < (\rho-|B|^2)(4x+\rho-|B|^2)$.
\end{enumerate}
\end{enumerate}
\end{proposition}

\begin{proof}
From Proposition~\ref{epsilon}, a Lie algebra $\frg$
underlying~(\ref{epsilonzero-red-new}) must be isomorphic to
$\frh_2$, $\frh_3$, $\frh_4$, $\frh_5$, $\frh_6$ or $\frh_8$.
Notice that the dimension of the center is~4 for $\frh_8$, 3~for
$\frh_6$, and~2 for the rest. The first Betti number is~5 for
$\frh_3$ and $\frh_8$, and~4 for $\frh_2,\frh_4,\frh_5$ and
$\frh_6$.

From Lemma~\ref{center}, $\frg$ is isomorphic to $\frh_6$ or
$\frh_8$ if and only if $|B|=\rho$ and $D=0$. Moreover, under
these conditions the first Betti number is~4 if $\rho=1$, and~5 if
$\rho=0$. So, (i.4) and (i.5) are proved.

Notice that $\frg$ has first Betti number equal to~5 if and only
if $B=\rho=0$ and $D\in\R$ in
equations~(\ref{epsilonzero-red-new}). Therefore, $\frg$ is
isomorphic to $\frh_3$ for $B=\rho=y=0$ and $x\not=0$, which
proves (i.2).

For the remaining cases $|B|\not=\rho$, $|B|=\rho$ and $y\not=0$,
or $|B|=\rho=1, y=0$ and $x\not=0$, the NLA $\frg$ has always
2-dimensional center by Lemma~\ref{center}, and its first Betti
number is equal to~4. Therefore, $\frg\cong\frh_2,\frh_4$ or
$\frh_5$. In order to decide which one is the corresponding Lie
algebra in terms of the coefficients $\rho$, $B$ and $D$, we
observe the following fact. Let $\alpha(\frg)$ be the number of
linearly independent elements $\tau$ in $\bigwedge ^2(\frg^*)$
such that $\tau\in d(\frg^*)$ and $\tau\wedge\tau=0$. It is
straightforward to check that $\alpha(\frh_k)$, for $k=2,4,5$,
equals the number of linearly independent exact 2-forms which are
decomposable, that is, $\alpha(\frh_2)=2$, $\alpha(\frh_4)=1$ and
$\alpha(\frh_5)=0$.

Let $\tau=\lambda\, d\omega^3 + \mu\, d\omega^{\bar{3}}$, where
$\lambda,\mu\in\C$, be any exact 2-form on $\frg$. Since $\tau$ is
real, $\mu=\bar{\lambda}$ and therefore
$$
\tau=\rho\lambda\,\omega^{12} +
(\lambda-\bar{\lambda})\omega^{1\bar{1}} +
B\lambda\,\omega^{1\bar{2}}
-\bar{B}\bar{\lambda}\,\omega^{2\bar{1}} +
(D\lambda-\bar{D}\bar{\lambda})\omega^{2\bar{2}} +
\rho\bar{\lambda}\,\omega^{\bar{1}\bar{2}}.
$$
A direct calculation shows that
$$
\tau\wedge\tau= 2\left( (\rho^2-|B|^2)|\lambda|^2 -(\lambda -
\bar{\lambda})(D\lambda-\bar{D}\bar{\lambda}) \right)\,
\omega^{12\bar{1}\bar{2}}.
$$
Thus, if we denote $p=\Re \lambda$ and $q=\Im \lambda$, then
$\tau\wedge\tau=0$ if and only if
\begin{equation}\label{second-deg}
(\rho-|B|^2)\, p^2 + 4ypq + (\rho-|B|^2 + 4x)\, q^2 = 0.
\end{equation}
If $|B|=\rho$ then (\ref{second-deg}) becomes $4q(yp+xq)=0$.
Therefore, $\tau_1=d(\Re \omega^3)$ is an exact 2-form on $\frg$
which is nonzero if $\rho=1$ or $y\not=0$, and it satisfies
$\tau_1\wedge\tau_1=0$. Moreover, when $\rho=1$, $y=0$ and
$x\not=0$, it follows from~(\ref{second-deg}) that $q=0$ and any
exact 2-form $\tau$ satisfying $\tau\wedge\tau=0$ must be a
multiple of $\tau_1$, thus $\alpha(\frg)=1$ and $\frg$ is
isomorphic to $\frh_4$, which proves (i.3). But when $y\not=0$,
the exact 2-form $\tau_2=-\frac{x}{y} d(\Re \omega^3) - d(\Im
\omega^3)$ satisfies $\tau_2\wedge\tau_2=0$. Since $\tau_1,\tau_2$
are linearly independent, we have that $\alpha(\frg)=2$ and
$\frg\cong\frh_2$, so (i.1) is proved. This completes the proof
of~(i).

To prove (ii), we consider (\ref{second-deg}) as a second-degree
equation in the variable $p$. Notice that the discriminant is
$\Delta=4q^2\left( 4y^2 - (\rho-|B|^2)(4x+\rho-|B|^2) \right)$,
and that $q\not=0$ because otherwise~(\ref{second-deg}) reduces to
$p=0$ and therefore $\lambda$ would be zero. Therefore, if $4y^2 >
(\rho-|B|^2)(4x+\rho-|B|^2)$ then $\Delta>0$ and for each
$q\not=0$, there exist two distinct solutions $p_1$ and $p_2$
of~(\ref{second-deg}). In this case we have $\alpha(\frg)=2$ and
therefore the underlying Lie algebra is isomorphic to $\frh_2$,
which proves (ii.1). A similar argument gives (ii.2) and (ii.3).
\end{proof}

\subsection{Equivalence of complex structures}\label{equival}

Let $\frg$ be a Lie algebra endowed with two complex structures
$J$ and $J'$. We recall that $J$ and $J'$ are said to be {\it
equivalent} if there is an automorphism $F\colon
\frg\longrightarrow \frg$ of the Lie algebra such that
$J'=F^{-1}\circ J\circ F$, that is, $F$ is a linear automorphism
such that $F^*\colon \frg^*\longrightarrow \frg^*$ commutes with
the Chevalley-Eilenberg differential~$d$ and $F$ commutes with the
complex structures $J$ and $J'$. The latter condition is
equivalent to say that~$F^*$, extended to the complexified
exterior algebra, preserves the bigraduations induced by $J$ and
$J'$.

It is clear that the nilpotency condition for a complex structure
is invariant under equivalence, that is, if $J'$ is equivalent to
$J$ then $J$ is nilpotent if and only if $J'$ is.

\begin{proposition}\label{equivspecial}
Any nilpotent, resp. nonnilpotent, complex structure on a
$6$-dimensional NLA is equivalent to a complex structure defined
by~$(\ref{nilpotent})$, resp. by~$(\ref{nonnilpotent})$, in the
sense of Proposition~$\ref{J-char}$.
\end{proposition}

\begin{proof}
Notice that if $\frg^{1,0}_J$ and $\frg^{1,0}_{J'}$ denote the
$(1,0)$-subspaces of $\gc^*$ associated to two complex
structures~$J$ and $J'$, then they are equivalent if and only if
there is a $\C$-linear isomorphism $F^*\colon \frg^{1,0}_J
\longrightarrow \frg^{1,0}_{J'}$ such that $d \circ F^*=F^* \circ
d$. Therefore, the result follows from Theorem~\ref{J-red}.
\end{proof}

Lemma~\ref{epsilonzero-new} states that any (non complex
parallelizable) complex structure on a $2$-step NLA with first
Betti number $\geq 4$ is equivalent to one defined
by~$(\ref{epsilonzero-red-new})$ in the sense of
Proposition~$\ref{J-char}$. Moreover:

\begin{corollary}\label{allequiv}
On the Lie algebras $\frh_6$ and $\frh_8$, any two complex
structures are equivalent.
\end{corollary}

\begin{proof}
From~(i.5) in Proposition~\ref{general} we have that any complex
structure on $\frh_8$ is equivalent to the one defined
by~(\ref{epsilonzero-red-new}) with $\rho=B=D=0$, and~(i.4) shows
that any complex structure on $\frh_6$ is equivalent to one
defined by~(\ref{epsilonzero-red-new}) with $\rho=|B|=1$ and
$D=0$. Since $|B|=1$ there exists a nonzero $\lambda$ satisfying
$\bar\lambda\, B=\lambda$, and the change of basis given by
$\omega'^1=\lambda \, \omega^1$, $\omega'^2=\bar\lambda \,
\omega^2$ and $\omega'^3=|\lambda|^2 \, \omega^3$ allows us to
consider the coefficient $B=1$.
\end{proof}

Let $J_0^+$ and $J_0^-$ be the abelian complex structures on the
Lie algebra $\frh_3$ defined by
$$
d\omega^1 = d\omega^2 = 0,\quad d\omega^3 = \omega^{1\bar{1}}
\pm\, \omega^{2\bar{2}}\, .
$$

\begin{corollary}\label{J0+-}
Any complex structure on $\frh_3$ is equivalent to $J_0^+$ or
$J_0^-$.
\end{corollary}

\begin{proof}
By Proposition~\ref{general}~(i.2) any complex structure on
$\frh_3$ is equivalent to one defined
by~(\ref{epsilonzero-red-new}) with $\rho=B=0$ and $D\in
\R-\{0\}$, and we can normalize $D$ to be $1$ or $-1$ depending on
the sign of $D$.
\end{proof}

Notice that the orientation induced by $J_0^+$ is opposite to the
one induced by the structure~$J_0^-$.

\section{Strong K\"ahler with torsion geometry in six dimensions}

Let $(J,g)$ be a Hermitian structure on a $2n$-dimensional
manifold $M$, that is, $J$ is a complex structure on $M$ which is
orthogonal relative to the Riemannian metric $g$. We denote by
$\Omega$ the {\it fundamental $2$-form} of $(J,g)$, which is
defined by $\Omega(X,Y)=g(JX,Y)$, for any differentiable vector
fields $X,Y$ on $M$.

It is well-known that the integrability of $J$ produces a
decomposition of the exterior differential $d$ of $M$ as
$d=\partial+\db$, where $\partial=\pi_{*+1,*}\circ d$ and $\db$ is
the conjugate of $\partial$. Since $d^2=0$, we have
$\partial^2=\db^2=0$ and $\partial \db=-\db \partial$.

A Hermitian structure $(J,g)$ is called {\it strong K\"ahler with
torsion} (SKT for short) if $\partial\Omega$ is a nonzero
$\db$-closed form. In this case, we shall refer to $g$ as an {\it
SKT metric}. Notice that a Hermitian structure $(J,g)$ is SKT if
and only if $Jd\Omega$ is nonzero and closed, because
$\partial\db$ acts as $\frac{1}{2}idJd$ on forms of bidegree
(1,1).

\medskip

Let $J$ be a complex structure on a Lie algebra $\frg$. An inner
product $g$ on $\frg$ such that
$g(J\,\cdot,J\,\cdot)=g(\cdot,\cdot)$ will be called {\it
$J$-Hermitian metric} on $\frg$, and we shall refer to the
associated $\Omega$ as the {\it fundamental form} of the {\it
Hermitian structure} $(J,g)$ on $\frg$.

Since $J$ is integrable on $\frg$, the extended
Chevalley-Eilenberg differential $d\colon \bigwedge^* \,\gc^*
\longrightarrow \bigwedge^{*+1} \,\gc^*$ also decomposes as
$d=\partial+\db$, where $\partial=\pi_{p+1,q}\circ d\colon$
$\bigwedge^{p,q}(\frg^*) \longrightarrow
\bigwedge^{p+1,q}(\frg^*)$ and $\db$ is the conjugate of
$\partial$. Any $J$-Hermitian metric $g$ on $\frg$ for which
$\partial\Omega$ is a nonzero $\db$-closed form will be called
{\it SKT metric} on $\frg$, and we shall refer to the pair $(J,g)$
as an {\it SKT structure} on $\frg$.

If the simply-connected nilpotent Lie group $G$ corresponding to
an NLA $\frg$ has a discrete subgroup~$\Gamma$ such that $M=\nilm$
is compact, then any Hermitian, resp. SKT, structure $(J,g)$ on
$\frg$ will pass to a Hermitian, resp. SKT, structure on the
nilmanifold $M$. Such a structure on $M$ will be also denoted by
$(J,g)$ and we shall refer to it as an {\it invariant} Hermitian,
resp. {\it invariant} SKT, structure on~$M$.

\medskip

Suppose that the NLA $\frg$ has dimension 6 and fix a basis
$\{\omega^j\}_{j=1}^3$ for $\frg^{1,0}$. Then, in terms of this
basis any $J$-Hermitian metric $g$ on $\frg$ is expressed as
\begin{equation}\label{metric}
\begin{array}{rcl}
g=&\!\!\! r \!\!\!&\!\!\omega^1\#\omega^{\bar{1}} +
s\,\omega^2\#\omega^{\bar{2}} + t\,\omega^3\#\omega^{\bar{3}} -
i\left( u\,\omega^1\#\omega^{\bar{2}}
-\bar{u}\,\omega^2\#\omega^{\bar{1}}\right. \\
&\!\!\!+\!\!\!&\! \left.
v\,\omega^2\#\omega^{\bar{3}}-\bar{v}\,\omega^3\#\omega^{\bar{2}}
+z\,\omega^1\#\omega^{\bar{3}}-\bar{z}\,\omega^3\#\omega^{\bar{1}}
\right),
\end{array}
\end{equation}
where $r,s,t\in\R$ and $u,v,z\in\C$ must satisfy restrictions that
guarantee that $g$ is positive definite, i.e. $g(Z,\bar{Z})>0$ for
any nonzero $Z\in (\frg^{1,0})^*$. Therefore, $r>0$, $s>0$, $t>0$,
$rs>|u|^2$, $st>|v|^2$, $rt>|z|^2$ and
$rst+2\,\Re(i\bar{u}\bar{v}z)> t|u|^2+r|v|^2+s|z|^2$.

The fundamental 2-form $\Omega\in\bigwedge^2\, \frg^*$ of the
Hermitian structure $(J,g)$ is then given by
\begin{equation}\label{2forma}
\Omega=i (r\,\omega^{1\bar{1}} + s\,\omega^{2\bar{2}} +
t\,\omega^{3\bar{3}})
+u\,\omega^{1\bar{2}}-\bar{u}\,\omega^{2\bar{1}} +
v\,\omega^{2\bar{3}}-\bar{v}\,\omega^{3\bar{2}}+
z\,\omega^{1\bar{3}}-\bar{z}\,\omega^{3\bar{1}}.
\end{equation}

The following result is proved by a direct calculation, so we omit
the proof.

\begin{lemma}\label{partialOmega}
Let $(J,g)$ be a Hermitian structure on a $6$-dimensional NLA
$\frg$, and $\Omega$ its fundamental form.
\begin{enumerate}
\item[{\rm (i)}] If $J$ is nonnilpotent, then in terms of the
basis $\{\omega^j\}_{j=1}^3$ of $\frg^{1,0}$
satisfying~$(\ref{nonnilpotent})$, the $(2,1)$-form
$\partial\Omega$ is given by
\begin{eqnarray*}
\begin{array}{rcl}
\partial\Omega=  &\!\!\! - \!\!\!&\!\!
(\bar{A}v+ibz)\omega^{12\bar{1}} - ibEv\,\omega^{12\bar{2}}
-(i\bar{A}t - u + E\bar{u})\omega^{13\bar{1}}\\[7pt]
&\!\!\! + \!\!\!&\!\! (is+bt)E\,\omega^{13\bar{2}} +
Ev\,\omega^{13\bar{3}} + (is-bt)\omega^{23\bar{1}}.
\end{array}
\end{eqnarray*}
\item[{\rm (ii)}] If $J$ is nilpotent, then in terms of the basis
$\{\omega^j\}_{j=1}^3$ of $\frg^{1,0}$
satisfying~$(\ref{nilpotent})$, the form $\partial\Omega$ is given
by
\begin{eqnarray*}
\begin{array}{rcl}
\partial\Omega=  &\!\!\! - \!\!\!&\!\!
\left( i\epsilon s + \rho \bar{z} + (1-\epsilon)\bar{A}v -
\bar{B}z \right)\omega^{12\bar{1}} - \left( \rho \bar{v} + \bar{C}
v - (1-\epsilon)\bar{D}z \right)\omega^{12\bar{2}} + i\rho t\,
\omega^{12\bar{3}}\\[7pt]
&\!\!\! + \!\!\!&\!\! \left( \epsilon \bar{v} -
i(1-\epsilon)\bar{A} t \right)\omega^{13\bar{1}}
-i\bar{C}t\,\omega^{13\bar{2}} - i\bar{B}t\,\omega^{23\bar{1}}
-i(1-\epsilon)\bar{D}t\,\omega^{23\bar{2}}.
\end{array}
\end{eqnarray*}
\end{enumerate}
\end{lemma}

\medskip

The theorem below is essentially given by Fino, Parton and Salamon
in~\cite[Theorems~1.2 and~3.2]{FPS}. Their proof involves a direct
but rather long calculation following a decision tree to eliminate
$B_{13}$, $B_{1\bar{3}}$, $B_{23}$, $B_{2\bar{3}}$ and the five
coefficients $A$'s in the general equations~(\ref{ecus}) under the
SKT hypothesis. We give a simple proof based on our previous study
of complex geometry developed in Section~2, together with the fact
that the SKT condition is satisfied up to equivalence of the
complex structure. Our proof also illustrates a general procedure
that is useful to investigate balanced and locally conformal
K\"ahler geometry, as it is shown in the next sections. Notice
that part~({\it a}) of the following theorem is a sligthly
stronger version of Theorem~1.2 in~\cite{FPS}.

\begin{lemma}\label{SKTinvariance}
Let $\frg$ be a Lie algebra endowed with a complex structure $J$
having compatible SKT metrics. Then, any other complex structure
$J'$ equivalent to $J$ posseses compatible SKT metrics.
\end{lemma}

\begin{proof}
Let $(J,g)$ be an SKT structure with fundamental form $\Omega$,
and $F\in {\rm Aut}(\frg)$ an automorphism such that $F\circ
J'=J\circ F$. Then, $g'=F^*g$ is a $J'$-Hermitian metric on $\frg$
with fundamental form $\Omega'=F^*\Omega$. Since $F^*$ commutes
with $d$ and preserves the bidegree, we get
$\db'\partial'\Omega'=F^*(\db\partial\Omega)$, where
$d=\partial'+\db'$ is the decomposition of $d$ with respect to
$J'$. Therefore, $\partial\Omega$ is a nonzero $\db$-closed form
if and only if $\partial'\Omega'$ is a nonzero $\db'$-closed form.
\end{proof}

\begin{theorem}\label{strongKT}
Let $\frg$ be a $6$-dimensional NLA admitting complex structures.
\begin{enumerate}
\item[{\rm ({\it a})}] A Hermitian structure $(J,g)$ on $\frg$ is
SKT if and only if the complex structure $J$ is equivalent to one
defined by~$(\ref{epsilonzero-red-new})$ with
\begin{equation}\label{SKTcondition}
\rho+|B|^2=2\,\Re(D).
\end{equation}
In particular, if $(J,g)$ is an SKT structure then $J$ is
nilpotent and any other $J$-Hermitian metric on $\frg$ is~SKT.
\item[{\rm ({\it b})}] There exists an SKT structure on $\frg$ if
and only if it is isomorphic to $\frh_{2}$, $\frh_{4}$, $\frh_{5}$
or $\frh_{8}$.
\end{enumerate}
\end{theorem}

\begin{proof}
To prove ({\it a}), we use Lemma~\ref{SKTinvariance} and
Proposition~\ref{equivspecial} to focus on the two special types
of complex structures defined by~(\ref{nonnilpotent})
and~(\ref{nilpotent}). If $J$ is a nonnilpotent complex structure
defined by~(\ref{nonnilpotent}), then it follows from
Lemma~\ref{partialOmega}~(i) that
$$
\db\partial\Omega = 2i(b^2 t\, \omega^{12\bar{1}\bar{2}} +s\,
\omega^{13\bar{1}\bar{3}})\not= 0,
$$
because $g$ is positive definite and in particular $s>0$. Thus,
$J$ must be necessarily nilpotent if it has a compatible SKT
metric, so it can be expressed by equations of the
form~(\ref{nilpotent}). Now, from Lemma~\ref{partialOmega}~(ii) we
have
$$
\db\partial\Omega = it \left(\rho^2+|B|^2+|C|^2 -
2(1-\epsilon)^2\, \Re(A\bar{D}) \right) \omega^{12\bar{1}\bar{2}}.
$$
If $\epsilon=1$ then we must have $\rho=B=C=0$ because $t>0$, so
in such case we can suppose $\epsilon=0$ after interchanging
$\omega^2$ with $\omega^3$. Also notice that a complex
parallelizable structure cannot satisfy the condition
$\db\partial\Omega=0$, unless the NLA $\frg$ be abelian, in which
case $\Omega$ would be closed. Therefore, we can apply
Proposition~\ref{epsilon} and Lemma~\ref{epsilonzero-new} to get
the equivalent condition~(\ref{SKTcondition}).

In order to prove ({\it b}), we first observe that
Proposition~\ref{general} implies that the possible candidates to
admit an SKT structure are $\frh_2,\ldots,\frh_6$ and $\frh_8$.
But, from~(i.2) and~(i.4) in Proposition~\ref{general} it follows
that there is no SKT structure neither on $\frh_3$ nor on
$\frh_6$, because~(\ref{SKTcondition}) is never satisfied. On the
other hand, (i.5) shows that any complex structure on $\frh_8$ has
compatible SKT metrics. Finally, the
condition~(\ref{SKTcondition}) for $\rho=1$ and $B=0$ in
equations~(\ref{epsilonzero-red-new}) reduces to
$D=\frac{1}{2}+iy$, so Proposition~\ref{general}~(ii) implies the
existence of SKT structures on $\frh_2$ for $|y|>
\frac{\sqrt{3}}{2}$, on $\frh_4$ for $y= \pm \frac{\sqrt{3}}{2}$,
and on $\frh_5$ for $|y|< \frac{\sqrt{3}}{2}$.
\end{proof}

\medskip

Next we describe (a parametrization of) the space of SKT
structures in dimension 6 up to equivalence of the underlying
complex structure. In view of~({\it a}) in the theorem above, we
restrict our attention to complex structures $J$ defined
by~(\ref{epsilonzero-red-new}) with $B=p+iq,D=x+iy \in \C$ and
$\rho=0,1$, and satisfying the SKT condition
$x=\frac{1}{2}(\rho+p^2+q^2)$. Let us fix $\rho=0$ or~$1$, which
is equivalent to require that $J$ be abelian or not. Then, the
complex structures having compatible SKT metric are parametrized
by points $(p,q,y)$ in the Euclidean space $\R^3$. Now, given an
NLA $\frg$ admitting SKT structures, we shall find which is the
region in the Euclidean space that parametrizes the space of
complex structures (up to equivalence) on $\frg$ satisfying the
SKT condition. For that, we make use of Proposition~\ref{general}
taking into account that $(\rho-|B|^2)(4x+\rho-|B|^2)=
4\rho-(\rho+p^2+q^2)^2$ under the SKT assumption:
\begin{enumerate}
\item[$\bullet$] First, let us suppose that $J$ is abelian, that
is, $\rho=0$. If $p=q=0$ then the corresponding Lie algebra is
$\frh_8$ for $y=0$, and $\frh_2$ for $y\not=0$. If $p^2+q^2\not=
0$ then $4y^2+(p^2+q^2)^2$ is strictly positive, so the
corresponding Lie algebra is again $\frh_2$. Therefore, the SKT
structures $(J,g)$ with abelian $J$ are parametrized by the points
in the Euclidean space $\R^3$, where the origin corresponds to SKT
structures on the Lie algebra $\frh_8$ and the points in
$\R^3-\{0\}$ to SKT structures on $\frh_2$.

\item[$\bullet$] Suppose now that $J$ is nilpotent but nonabelian,
i.e. $\rho=1$. If $p^2+q^2=1$ then the corresponding Lie algebra
is $\frh_4$ for $y=0$, and $\frh_2$ for $y\not=0$. If
$p^2+q^2\not= 1$ then the equation $4y^2-4+(1+p^2+q^2)^2=0$
represents an ovaloid of revolution generated by rotating the
curve illustrated in the Figure~1 about the $y$-axis. Therefore,
the SKT structures $(J,g)$ with nonabelian $J$ are parametrized by
the points in the Euclidean 3-space, where the region outside the
ovaloid corresponds to SKT structures on the Lie algebra $\frh_2$,
the points on the ovaloid correspond to SKT structures on $\frh_4$
and the region inside the ovaloid to SKT structures on $\frh_5$.

\bigskip

\setlength{\unitlength}{2.5pt}
\begin{picture}(120,50)(-82,-18)
\multiput(-17.5,1)(0,2){5}{\circle*{.2}}
\multiput(-15,0)(0,2){7}{\circle*{.2}}
\multiput(-12.5,1)(0,2){8}{\circle*{.2}}
\multiput(-10,0)(0,2){9}{\circle*{.2}}
\multiput(-7.5,1)(0,2){9}{\circle*{.2}}
\multiput(-5,0)(0,2){9}{\circle*{.2}}
\multiput(-2.5,1)(0,2){9}{\circle*{.2}}
\multiput(0,0)(0,2){10}{\circle*{.2}}
\multiput(2.5,1)(0,2){9}{\circle*{.2}}
\multiput(5,0)(0,2){9}{\circle*{.2}}
\multiput(7.5,1)(0,2){9}{\circle*{.2}}
\multiput(10,0)(0,2){9}{\circle*{.2}}
\multiput(12.5,1)(0,2){8}{\circle*{.2}}
\multiput(15,0)(0,2){7}{\circle*{.2}}
\multiput(17.5,1)(0,2){6}{\circle*{.2}}
\multiput(20,0)(0,2){3}{\circle*{.2}}
\multiput(-17.5,1)(0,-2){6}{\circle*{.2}}
\multiput(-15,0)(0,-2){7}{\circle*{.2}}
\multiput(-12.5,1)(0,-2){9}{\circle*{.2}}
\multiput(-10,0)(0,-2){9}{\circle*{.2}}
\multiput(-7.5,1)(0,-2){10}{\circle*{.2}}
\multiput(-5,0)(0,-2){9}{\circle*{.2}}
\multiput(-2.5,1)(0,-2){10}{\circle*{.2}}
\multiput(0,0)(0,-2){10}{\circle*{.2}}
\multiput(2.5,1)(0,-2){10}{\circle*{.2}}
\multiput(5,0)(0,-2){9}{\circle*{.2}}
\multiput(7.5,1)(0,-2){10}{\circle*{.2}}
\multiput(10,0)(0,-2){9}{\circle*{.2}}
\multiput(12.5,1)(0,-2){9}{\circle*{.2}}
\multiput(15,0)(0,-2){7}{\circle*{.2}}
\multiput(17.5,1)(0,-2){7}{\circle*{.2}}
\multiput(20,0)(0,-2){3}{\circle*{.2}}
\put(0,0){\vector(1,0){30}} \put(0,0){\vector(0,1){30}}
\put(30,-3){$p=\Re B$} \put(2,28){$y=\Im D$}
\put(-25,20){$\frh_2$} \put(15,15){$\frh_4$}
\put(-11.65,-1){$\frh_5$}
\linethickness{1pt}\qbezier(20,0)(20,18)(0,18)
\qbezier(-20,0)(-20,18)(0,18) \qbezier(20,0)(20,-18)(0,-18)
\qbezier(-20,0)(-20,-18)(0,-18)
\end{picture}

\bigskip

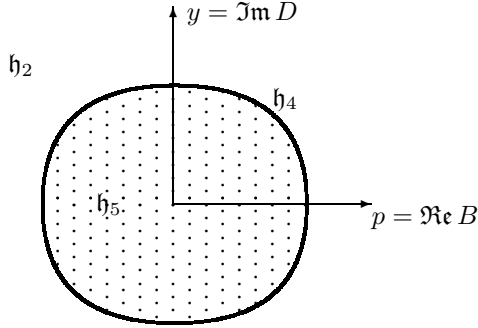
\centerline{{\bf Figure 1}: SKT structures in the nonabelian case}

\bigskip

\end{enumerate}

The Lie algebras $\frh_{2}$, $\frh_{4}$, $\frh_{5}$ and $\frh_{8}$
posses abelian complex structures. The following result is a
direct consequence of our study above.

\begin{corollary}\label{abelianSKT}
Let $(J,g)$ be an SKT structure on a $6$-dimensional NLA $\frg$.
If $J$ is abelian then $\frg$ is isomorphic to $\frh_{2}$ or
$\frh_{8}$. Therefore, none of the abelian complex structures on
$\frh_{4}$ or $\frh_{5}$ admits compatible SKT metric.
\end{corollary}

\medskip

Next we prove that the existence of an SKT structure on a compact
nilmanifold $\nilm$ implies the existence of an invariant one. The
proof is based on the ``symmetrization'' process, and we follow
ideas of \cite{Belg,FG}.

Let $M=\nilm$ be a compact nilmanifold and $\nu=d\tau$ a volume
element on $M$ induced by a bi-invariant one on the Lie group
$G$~\cite{Milnor}. After rescaling, we can suppose that $M$ has
volume equal to 1.

Given any covariant $k$-tensor field $T\colon {\frak X}(M)\times
\cdots \times {\frak X}(M) \longrightarrow {\cal C}^{\infty}(M)$
on the nilmanifold $M$, we define a covariant $k$-tensor $T_{\nu}
\colon \frg\times \cdots \times \frg \longrightarrow \R$ on the
NLA $\frg$ of $G$ by
$$
T_{\nu}(X_1,\ldots,X_k)=\int_{p\in M} T_p
(X_1\!\mid_p,\ldots,X_k\!\mid_p) \, \nu\ , \quad\quad \mbox{for}\
X_1,\ldots,X_k\in \frg,
$$
where $X_j\!\mid_p$ is the value at the point $p\in M$ of the
projection on $M$ of the left-invariant vector field $X_j$ on the
Lie group $G$.

Obviously, $T_{\nu}=T$ for any tensor field $T$ coming from a
left-invariant one. In~\cite{Belg} it is proved that if $T=\alpha
\in {\cal A}^k(M)$ is a $k$-form on $M$, then
$(d\alpha)_{\nu}=d\alpha_{\nu}$. A simple calculation shows that
$(\alpha_\nu\wedge\beta)_\nu=\alpha_\nu\wedge\beta_\nu$, for any
$\alpha \in {\cal A}^k(M)$ and $\beta \in {\cal A}^l(M)$.

\begin{remark}\label{nomizu}
{\rm The symmetrization process defines a linear map
$\tilde{\nu}\colon {\cal A}^k(M) \longrightarrow
\bigwedge\phantom{\!\!}^k \,\frg^*$, given by
$\tilde{\nu}(\alpha)=\alpha_\nu$ for any $\alpha \in {\cal
A}^k(M)$, which commutes with the differentials. Moreover, it
follows from Nomizu theorem~\cite{No} that $\tilde{\nu}$ induces
an isomorphism $H^k(\tilde{\nu})\colon H^k(M) \longrightarrow
H^k(\frg)$ between the $k$th cohomology groups, for each~$k$. In
particular, any closed $k$-form $\alpha$ on $M$ is cohomologous to
the invariant $k$-form $\alpha_\nu$ obtained by the symmetrization
process.}
\end{remark}

Let us suppose now that the nilmanifold $M=\nilm$ is equipped with
an invariant complex structure~$J$. If $g$ is a $J$-Hermitian
metric on $M$ and $\Omega$ denotes its fundamental form, then
$g_{\nu}$ is a $J$-Hermitian metric on the NLA $\frg$ with
fundamental form $\Omega_{\nu}$.

\begin{proposition}\label{symmetrization-SKT}
Let $(M=\nilm,J)$ be a compact complex nilmanifold with invariant
$J$, and suppose that the NLA $\frg$ of $G$ is not abelian. If
$(J,g)$ is an SKT structure on $M$ then $(J,g_{\nu})$ is an SKT
structure on $\frg$.
\end{proposition}

\begin{proof}
The fundamental form $\Omega$ of $(J,g)$ satisfies $dJd\Omega=0$,
but $d\Omega\not=0$. As it is observed in~\cite{FG}, since $J$ is
invariant, we have that $(J\alpha)_{\nu}=J\, \alpha_{\nu}$ for any
$\alpha\in{\cal A}^k(M)$. Therefore,
$dJd\Omega_{\nu}=dJ(d\Omega)_{\nu}=d(Jd\Omega)_{\nu}=(dJd\Omega)_{\nu}=0$.
Moreover, since $\frg$ is not abelian, it follows from~\cite{BG}
that $d\Omega_{\nu}\not=0$ because $M$ has no K\"ahler metric.
\end{proof}

Notice that the symmetrization of SKT structures on a torus gives
rise to invariant K\"ahler metrics.

\begin{corollary}\label{SKT-general}
A non-toral compact nilmanifold $M=\nilm$ of dimension~$6$ admits
an SKT metric compatible with an invariant complex structure if
and only if the Lie algebra of $G$ is isomorphic to $\frh_2,
\frh_4,\frh_5$ or $\frh_8$.
\end{corollary}

The result follows directly from Theorem~\ref{strongKT}~({\it b})
and Proposition~\ref{symmetrization-SKT}. In particular, the first
Betti number of $M$ must be $\geq 4$ in order to admit an SKT
structure $(J,g)$ with invariant $J$.

\medskip

Finally, Corollary~\ref{abelianSKT} implies that $M=\nilm$ has SKT
structures whose underlying complex structure is abelian if and
only if the Lie algebra of $G$ is isomorphic to $\frh_2$ or
$\frh_8$.

\section{Balanced metrics on six dimensional nilmanifolds}

Let $(J,g)$ be a Hermitian structure on a $2n$-dimensional
manifold $M$ with fundamental form $\Omega$. According
to~\cite{GH}, the $W_1$ and $W_2$ components in the well-known
Gray-Hervella decomposition of $\nabla\Omega$ are identically
zero, that is, $\nabla\Omega\in W_3\oplus W_4$. In this section we
are interested in Hermitian structures satisfying $\nabla\Omega\in
W_3$.

Let $\theta$ be the {\it Lee form} associated to the Hermitian
structure $(J,g)$, that is, $\theta=\frac{1}{1-n} J\delta\Omega$,
where $\delta$ denotes the formal adjoint of $d$ with respect to
the metric $g$. If $\theta$ vanishes identically then the
Hermitian structure is called {\it balanced} and we shall say that
$g$ is a {\it balanced} metric on~$M$. Such structures correspond
to the Gray-Hervella class ${\cal W}_3$~\cite{GH}, and they are
also known in the literature as {\it cosymplectic} or {\it
semi-K\"ahler}.

A large class of balanced structures is provided by the compact
{\it complex parallelizable} manifolds, that is, compact complex
manifolds $M$ for which the holomophic bundle $T^{1,0}M$ is
trivial. Wang~\cite{Wang} proved that $M$ is a compact quotient
$\Gamma\backslash G$ of a simply connected {\it complex} Lie group
$G$ by a discrete subgroup~$\Gamma$. Therefore, $G$ is
unimodular~\cite{Milnor}, so any Hermitian left-invariant metric
on the complex Lie group $G$ is balanced by~\cite[Theorem 2.2]{AG}
and it descends to $M$.

Alexandrov and Ivanov prove in~\cite[Remark 1]{AI} (see
also~\cite[Proposition 1.4]{FPS}) that the balanced condition is
complementary to the SKT condition in dimension~$\geq 6$. As a
consequence we have:

\begin{proposition}\label{cpar}
A compact complex parallelizable manifold (not a torus) of
dimension~$\geq 6$  has no compatible SKT metrics.
\end{proposition}

\begin{proof}
Let $M=\Gamma\backslash G$ be a compact complex parallelizable
manifold and denote by $J$ its natural complex structure. Any
Hermitian left-invariant metric on $G$ does not satisfy the SKT
condition, because it is balanced. So there are no left-invariant
SKT metrics on $G$ compatible with $J$.

Moreover, since $G$ is unimodular there exists a bi-invariant
volume element, and the symmetrization of an SKT metric on $M$
would be a left-invariant SKT metric on $G$. In fact, the proof of
Proposition~\ref{symmetrization-SKT} is valid in this context,
except that we use Theorem 2.1 in~\cite{AG}, which states that if
$G$ is not abelian then there are no left-invariant K\"ahler
metrics on $G$ compatible with $J$, in order to ensure that the
symmetrization of the fundamental form is not closed.
\end{proof}

Let $\frg$ be a Lie algebra of dimension~$2n$ endowed with a
Hermitian structure $(J,g)$, in the sense of Section~3, with Lee
form $\theta\in\frg^*$. We say that $(J,g)$ is a {\it balanced}
structure, or that $g$ is a {\it balanced} metric, on $\frg$ if
$\theta=0$.

Fixed a complex structure $J$ on $\frg$, we denote by ${\cal
M}_3(\frg,J)$ the set of all balanced $J$-Hermitian metrics~$g$ on
$\frg$.

\begin{lemma}\label{Binvariance}
If $J'$ is a complex structure on $\frg$ equivalent to $J$, then
the metrics in ${\cal M}_3(\frg,J')$ are in one-to-one
correspondence with the metrics in ${\cal M}_3(\frg,J)$.
\end{lemma}

\begin{proof}
Let $F\in {\rm Aut}(\frg)$ be an automorphism of the Lie algebra
such that $F\circ J'=J\circ F$. Given $g\in {\cal M}_3(\frg,J)$
with fundamental form $\Omega$, let us consider the $J'$-Hermitian
metric $g'=F^*g$, whose fundamental form is $\Omega'=F^*\Omega$.
If we denote by $\delta'$ the adjoint of $d$ with respect to the
metric $g'$ then $\delta' F^*=F^* \delta$, which implies that the
Lee form $\theta'$ of the Hermitian structure $(J',g')$ is given
by $\theta'=F^*\theta$. Therefore, $(J,g)$ is balanced if and only
if $(J',g')$ is.
\end{proof}

When $\frg$ is 6-dimensional, $2*\Omega=\Omega\wedge\Omega$, where
$*$ denotes the Hodge star with respect to~$g$. So the Lee form
vanishes if and only if $\Omega^2$ is closed. But,
$d\Omega^2=2\Omega\wedge d\Omega$ is a real 5-form which
decomposes as a sum of forms of types $(3,2)$ and $(2,3)$, because
the bidegree of $\Omega$ is equal to $(1,1)$. Thus, $\Omega^2$ is
closed if and only if $(d\Omega^2)^{3,2}=2\Omega\wedge
(d\Omega)^{2,1}=0$. Therefore, a Hermitian structure is balanced
if and only if $\partial\Omega\wedge\Omega=0$.

Fixed a complex structure $J$ on an NLA $\frg$ of dimension~6, the
set ${\cal M}_3(\frg,J)$ is then given by
$${\cal M}_3(\frg,J)=\{ J\mbox{-Hermitian metrics } g\mid
\partial\Omega_g\wedge\Omega_g=0\},$$
where $\Omega_g$ is the fundamental form associated to $g$. Our
first goal is to prove that ${\cal M}_3(\frg,J)\not= \emptyset$
only for a Lie algebra $\frg$ isomorphic to
$\frh_{1},\ldots,\frh_6$ or $\frh_{19}^-$.

\begin{proposition}\label{balanced}
Let $(J,g)$ be a Hermitian structure on a $6$-dimensional NLA
$\frg$.
\begin{enumerate}
\item[{\rm ({\it a})}] If $J$ is nonnilpotent, then $(J,g)$ is
balanced if and only if the complex structure $J$ is equivalent to
one defined by~$(\ref{nonnilpotent})$ and the metric coefficients
of $g$ in~$(\ref{metric})$ satisfy
\begin{equation}\label{Bnonnilp-condition}
z=\frac{-iuv}{s} \quad\quad \mbox{and} \quad\quad
As+b\bar{E}u+b\bar{u}=0.
\end{equation}
\item[{\rm ({\it b})}] If $J$ is nilpotent but not complex
parallelizable, then $(J,g)$ is balanced if and only if $J$ is
equivalent to one defined by~$(\ref{epsilonzero-red-new})$ and the
metric coefficients of $g$ in~$(\ref{metric})$ satisfy
\begin{equation}\label{Bnilp-condition}
st-|v|^2 + D(rt-|z|^2) = B(it\bar{u} - v\bar{z}).
\end{equation}
\item[{\rm ({\it c})}] If $J$ is a complex parallelizable
structure, then any $J$-Hermitian metric is balanced.
\end{enumerate}
\end{proposition}

\begin{proof}
Suppose first that $J$ is nonnilpotent. From
Lemma~\ref{Binvariance} and Proposition~\ref{equivspecial}, we can
restrict our attention to fundamental 2-forms $\Omega$ given
by~(\ref{2forma}) in terms of a basis $\{\omega^j\}_{j=1}^3$
satisfying~(\ref{nonnilpotent}). Using
Lemma~\ref{partialOmega}~(i), a direct calculation shows that
$$
\partial\Omega\wedge\Omega =  \left( \bar{A}(st-|v|^2) +
b(tu-i\bar{v}z) + bE(t\bar{u}+iv\bar{z}) \right)
\omega^{123\bar{1}\bar{2}} + (uv-isz) \omega^{123\bar{1}\bar{3}}.
$$
Therefore, a metric $g$ given by~(\ref{metric}) is balanced if and
only if $z=-iuv/s$ and
$$
0=\bar{A}(st-|v|^2) + b(tu-i\bar{v}z) +
bE(t\bar{u}+iv\bar{z})={st-|v|^2\over s} (\bar{A}s + bu +
bE\bar{u}).
$$
Since $g$ is positive definite, the latter condition is equivalent
to
$$
As+b\bar{E}u+b\bar{u}=0
$$
because $s$ and $b$ are real numbers, so part ({\it a}) of the
proposition is proved.

To prove ({\it b}) we can focus, again by Lemma~\ref{Binvariance}
and Proposition~\ref{equivspecial}, on nilpotent complex
structures~$J$ defined by equations of the form~(\ref{nilpotent}).
For any $\Omega$ given by~(\ref{2forma}), from
Lemma~\ref{partialOmega}~(ii) we get by a simple calculation that
$$
\begin{array}{rcl}
\partial\Omega\wedge\Omega \! &\!\!\! = \!\!\!&\!\!
\left( (1\!-\!\epsilon)\bar{A}(st\!-\!|v|^2) +
\bar{B}(itu+\bar{v}z) - \bar{C}(it\bar{u}-v\bar{z}) +
(1\!-\!\epsilon)\bar{D}(rt\!-\!|z|^2)
\right) \omega^{123\bar{1}\bar{2}}\\[10pt]
&\!\!\! \phantom{-} \!\!\!&\!\! - \epsilon (st-|v|^2)
\omega^{123\bar{1}\bar{3}}.
\end{array}
$$
Since $g$ is positive definite, the coefficient of
$\omega^{123\bar{1}\bar{3}}$ vanishes if and only if $\epsilon=0$.
Thus, if $J$ is not complex parallelizable, then
Proposition~\ref{epsilon} and Lemma~\ref{epsilonzero-new} imply
that $J$ is equivalent to one defined
by~(\ref{epsilonzero-red-new}), and so the form
$\partial\Omega\wedge\Omega$ is zero if and only
if~(\ref{Bnilp-condition}) holds.

Finally, if $\epsilon=A=B=C=D=0$ then $\partial\Omega\wedge\Omega$
vanishes identically, so ({\it c}) is clear. It also follows
directly from~\cite{AG}.
\end{proof}

\begin{theorem}\label{balanced-clasif}
Let $\frg$ be an NLA of dimension~6. Then, there exists a balanced
structure on $\frg$ if and only if it is isomorphic to $\frh_{k}$,
for $1\leq k\leq 6$, or $\frh_{19}^-$. Moreover:
\begin{enumerate}
\item[{\rm ({\it a})}] Any complex structure on $\frh_6$ and
$\frh_{19}^-$ has compatible metrics which are balanced.
\item[{\rm ({\it b})}] A complex structure on $\frh_3$ has
balanced compatible metrics iff it is equivalent to~$J_0^-$.
\item[{\rm ({\it c})}] On the Lie algebras $\frh_2$, $\frh_4$ and
$\frh_5$ there exist complex structures having balanced compatible
metrics, but there also exist complex structures not admitting
such metrics.
\end{enumerate}
\end{theorem}

\begin{proof}
If there exists a balanced structure $(J,g)$ on $\frg$ such that
$J$ is nonnilpotent, then it follows
from~(\ref{Bnonnilp-condition}) by complex conjugation that
$$
\bar{A}s+bu+bE\bar{u}=0.
$$
On the other hand, if we multiply the second equation
in~(\ref{Bnonnilp-condition}) by $E$, then taking into account
that $|E|=1$ we get
$$
AEs+bu+bE\bar{u}=0.
$$
Therefore, $s(\bar{A}-AE)=0$, that is, $\bar{A}=AE$ because $g$ is
positive definite. Now Proposition~\ref{h19-h26} implies that
$\frg$ cannot be isomorphic to $\frh_{26}^+$.

Now suppose that $\frg$ has a balanced structure $(J,g)$ such that
$J$ is nilpotent. Propositions~\ref{epsilon} and~\ref{balanced}
imply that, up to isomorphism, the possible candidates for such a
Lie algebra are $\frh_1,\ldots,\frh_6$ and $\frh_8$. But the Lie
algebra $\frh_8$ is excluded by Proposition~\ref{general}~(i.5),
because~(\ref{Bnilp-condition}) reduces to $st-|v|^2 = 0$ for
$B=D=0$, which contradicts that $g$ is positive definite.
Therefore, $\frg$ cannot be isomorphic to $\frh_{k}$ for~$7\leq
k\leq 16$.

Notice that for the ``canonical'' metric $g$ given by $r=s=t=1$
and $u=v=z=0$, the balanced condition~(\ref{Bnilp-condition})
reduces to $D=-1$. From Proposition~\ref{general} it follows that
there is a balanced structure on $\frh_2$ for $|B|<1=\rho$, on
$\frh_4$ for $\rho=|B|=1$, on $\frh_5$ for $|B|>1=\rho$ and on
$\frh_3$ for the complex structure $J_0^-$, i.e. for $\rho=B=0$.

To complete the proof, it remains to show that any complex
structure on $\frh_6$ and $\frh_{19}^-$ has a compatible balanced
metric, and that there exists a complex structure on each one of
the Lie algebras $\frh_2$, $\frh_3$, $\frh_4$ and $\frh_5$
admitting no compatible balanced metric.

Let $g_u$ be the metric defined by $r=1+|u|^2$, $s=t=1$ and
$v=z=0$. If $u=-\bar{A}/(2b)$ then we have a metric $g_u$ on
$\frh_{19}^-$ compatible with the complex structure $J$ defined
by~(\ref{nonnilpotent}) with $E=\bar{A}/A$, according to
Proposition~\ref{h19-h26}. Since $g_u$
satisfies~(\ref{Bnonnilp-condition}), from Lemma~\ref{Binvariance}
we conclude that any other complex structure on $\frh_{19}^-$ has
a balanced compatible metric. Moreover, if $u=i$ then the metric
$g_u$ on $\frh_6$ is $J$-Hermitian for the complex structure $J$
defined by~(\ref{epsilonzero-red-new}) with $\rho=B=1$ and $D=0$,
according to Proposition~\ref{general}~(i.4), and it is clear
that~(\ref{Bnilp-condition}) holds. From Corollary~\ref{allequiv}
and Lemma~\ref{Binvariance} it follows that any other complex
structure on $\frh_6$ possesses a balanced compatible metric.

On the other hand, for the complex structure $J_0^+$ on $\frh_3$
given in Corollary~\ref{J0+-} the balanced
condition~(\ref{Bnilp-condition}) reduces to $st-|v|^2 + rt-|z|^2
=0$, so $g$ cannot be positive definite. Therefore, there is no
balanced compatible metric.

Finally, if $\rho=1$ and $B=x=0$ in Proposition~\ref{general}
then~(\ref{Bnilp-condition}) reduces to $st-|v|^2 + iy(rt-|z|^2)
=0$, so the metric cannot be positive definite, and depending on
the value of $y$ we get complex structures on the Lie algebras
$\frh_2$, $\frh_4$ and $\frh_5$ having no balanced compatible
metric.
\end{proof}

\begin{remark}
{\rm In~\cite{FPS} it is shown that the metric given by
$r=s=t=1/2$, $u=v=z=0$ is balanced with respect to one particular
complex structure on $\frh_k$, for $1\leq k\leq 6$.}
\end{remark}

By the symmetrization process, Fino and Grantcharov prove
in~\cite{FG} that if $J$ is an invariant complex structure on a
compact nilmanifold $M=\nilm$ admitting a balanced metric~$g$,
then there is a balanced structure $(J,\hat{g})$ on the Lie
algebra $\frg$ of $G$. Therefore:

\begin{corollary}\label{B-Jnn-nilmanifold}
A compact nilmanifold $M=\nilm$ of dimension~$6$ admits a balanced
metric compatible with an invariant complex structure if and only
if the Lie algebra of $G$ is isomorphic to $\frh_{19}^-$ or
$\frh_{k}$ for some $1\leq k\leq 6$.
\end{corollary}

Since $\frh_3$ is the Lie algebra underlying the compact
nilmanifold $N(2,1)\times S^1$, where $N(2,1)$ is a quotient of
the 5-dimensional generalized Heisenberg group $H(2,1)$, we have
that an invariant complex structure~$J$ on $N(2,1)\times S^1$ has
compatible balanced metrics if and only if $J$ is equivalent to
$J_0^-$.

\bigskip

Let $(J,g)$ be a Hermitian structure on a manifold $M$ and denote
by $\nabla^B$ its {\it Bismut connection} , i.e. the unique
connection for which $g$ and $J$ are parallel and the torsion
$T^B$ is given by $g(X,T^B(Y,Z))=Jd\Omega(X,Y,Z)$. Gutowski,
Ivanov and Papadopoulos pose in~\cite{GIP} the following version
of the Calabi conjecture in the non-K\"ahler setting: on {\it any}
$2n$-dimensional compact complex manifold with vanishing first
Chern class there exists a Hermitian structure with restricted
holonomy of the Bismut connection contained in ${\rm SU}(n)$. They
prove that this property holds for Moishezon manifolds, for
compact complex manifold with vanishing first Chern class which
are cohomologically K\"ahler and for connected sums of $k \geq 2$
copies of $S^3 \times S^3$.

Now, let $M=\nilm$ be a compact nilmanifold of dimension~6
equipped with an invariant complex structure $J$. It follows
from~(\ref{ecus}) that $\omega^{123}$ is a holomorphic
non-vanishing (3,0)-form. Therefore, Theorem~4.1 in~\cite{FG}
asserts that if the Ricci tensor of the Bismut connection of some
$J$-Hermitian metric $g$ on $M$ vanishes, then there is a metric
$\widetilde{g}$ conformal to $g$ such that $(J,\widetilde{g})$ is
a balanced structure on $M$, so there is an invariant balanced
structure $(J,\hat{g})$ on $M$ by the symmetrization process.
Conversely, given an invariant balanced Hermitian structure on $M$
there is a conformal change of metric such that the Ricci tensor
of the Bismut connection of the resulting metric vanishes
(see~\cite[Proposition 6.1]{FPS}).

Using this result, Fino and Grantcharov provide counter-examples
to the above conjecture by showing that there exist invariant
complex structures on the Iwasawa manifold which do not admit
compatible invariant balanced metrics. In the following result we
describe the general situation for 6-dimensional nilmanifolds.

\begin{corollary}\label{gen-Bismut}
Let $M=\nilm$ be a $6$-dimensional compact nilmanifold, and $\frg$
the Lie algebra of $G$. Then:
\begin{enumerate}
\item[{\rm ({\it a})}] If $\frg$ is isomorphic to $\frh_6$ or
$\frh_{19}^-$, then any invariant complex structure on $M$ has a
Hermitian structure with restricted holonomy of the Bismut
connection contained in ${\rm SU}(3)$. \item[{\rm ({\it b})}] If
$\frg$ is isomorphic to $\frh_2$, $\frh_3$, $\frh_4$ or $\frh_5$,
then there are invariant complex structures on $M$ having a
Hermitian structure with restricted holonomy of the Bismut
connection contained in ${\rm SU}(3)$, but there also exist
invariant complex structures on $M$ for which the restricted
holonomy of the Bismut connection of any Hermitian metric is not
contained in ${\rm SU}(3)$. \item[{\rm ({\it c})}] If $\frg$ is
isomorphic to $\frh_{26}^+$ or $\frh_{k}$ for some $7\leq k\leq
16$, then any invariant complex structure on $M$ does not posses
Hermitian structures with restricted holonomy of the Bismut
connection contained in ${\rm SU}(3)$.
\end{enumerate}
\end{corollary}

Observe that on the compact nilmanifold $N(2,1)\times S^1$, an
invariant complex structure $J$ has a Hermitian structure for
which the holonomy of its Bismut connection reduces to {\rm SU(3)}
if and only if $J$ is equivalent to $J_0^-$.

\medskip

We finish this section with some remarks about the stability of
the balanced condition under small deformations of the complex
structure. As a consequence of Theorem~\ref{balanced-clasif}, the
nilmanifolds corresponding to $\frh_6$ and $\frh_{19}^-$ are
stable in the sense that given a balanced structure $(J_0,g_0)$,
with $J_0$ invariant, then along any deformation $J_\alpha$ of
$J_0$ consisting of invariant complex structures, there always
exists a balanced $J_\alpha$-Hermitian metric $g_\alpha$ for each
value of $\alpha$.

However, it is shown in~\cite[Theorem 4.2]{FG} that such stability
does not hold for the Iwasawa manifold. Next we show that for the
compact nilmanifold $\Gamma\backslash (H^3\times H^3)$, where
$H^3$ is the Heisenberg group, the balanced condition is not
stable under small deformations.

Since the Lie algebra of $H^3\times H^3$ is $\frh_2$, we consider
the equations~(\ref{epsilonzero-red-new}) with $\rho=1$ and $B=0$,
and denote by $J_{x,y}$ the complex structure on $\frh_2$ given by
$$d \omega^1=d\omega^2=0,\quad d\omega^3=\omega^{12} + \omega^{1\bar{1}} +
(x+iy)\,\omega^{2\bar{2}},$$ where $x,y\in \R$ satisfy $4y^2-4x-1
> 0$, according to Proposition~\ref{general}. Observe that the balanced
condition~(\ref{Bnilp-condition}) is satisfied for $J_{x,y}$ if
and only if $x<-\frac{1}{4}$ and $y=0$. In particular, for $x=-1$
the family of complex structures $J_{-1,y}$ with $y\in\R$ only
admits balanced metrics for~$y=0$.

\begin{corollary}\label{nonstableh2}
The property ``vanishing Ricci tensor for the Bismut connection''
is not stable under small deformations on the nilmanifold
$\Gamma\backslash (H^3\times H^3)$.
\end{corollary}

\section{Locally conformal K\"ahler geometry}

A Hermitian structure $(J,g)$ on a $2n$-dimensional manifold $M$
is called {\it locally conformal K\"ahler} ({\it LCK} for short)
if $M$ has an open cover $\{U_i\}$ and a family $\{f_i\}$ of
differentiable functions $f_i\colon U_i\longrightarrow \R$ such
that each local metric $g_i=\exp f_i\, g_{|_{U_i}}$ is K\"ahler.
If $\Omega$ denotes the fundamental form, then the Hermitian
structure is LCK if and only if $d\Omega=\theta\wedge\Omega$ with
closed Lee form $\theta$. Notice that the class of LCK structures
corresponds to the Gray-Hervella class ${\cal W}_4$~\cite{GH}.

An interesting special class of locally conformal K\"ahler metrics
is the one consisting of those for which~$\theta$ is a nowhere
vanishing parallel form. Such Hermitian metrics are called {\it
generalized Hopf} metrics, and a Hermitian manifold endowed with
such a metric is also known as a {\it Vaisman} manifold
(see~\cite{DO,V2,V3}).

Let $(J,g)$ be a Hermitian structure on a Lie algebra $\frg$, with
fundamental form $\Omega\in\bigwedge^2\, \frg^*$ and Lee form
$\theta\in\frg^*$. We say that $(J,g)$ is a {\it LCK} structure,
or that $g$ is a {\it LCK} metric, if $d\Omega=\theta\wedge\Omega$
with closed Lee form $\theta$.

Fixed a complex structure $J$ on $\frg$, we denote by ${\cal
M}_4(\frg,J)$ the set of LCK $J$-Hermitian metrics on $\frg$.

\begin{lemma}\label{LCKinvariance}
If $J'$ is a complex structure on $\frg$ equivalent to $J$, then
the metrics in ${\cal M}_4(\frg,J')$ are in one-to-one
correspondence with the metrics in ${\cal M}_4(\frg,J)$.
\end{lemma}

\begin{proof}
Following the proof of Lemma~\ref{Binvariance}, we have
$F^*(d\Omega-\theta\wedge\Omega)= d F^*\Omega-F^*\theta\wedge
F^*\Omega= d\Omega'-\theta'\wedge\Omega'$.
\end{proof}

Since $\theta$ and $\Omega$ are real forms, taking into account
their bidegrees we have that in dimension $\geq 6$ a~Hermitian
structure is LCK if and only if $\partial\Omega
=\theta^{1,0}\wedge\Omega$. Therefore, if $\dim \frg \geq 6$ then
$${\cal M}_4(\frg,J)=\{ J\mbox{-Hermitian metrics } g\mid
\partial\Omega_g - (\theta_g)^{1,0}\wedge\Omega_g = 0\}.
$$

\begin{theorem}\label{LCKgen}
A $6$-dimensional NLA $\frg$ admits an LCK structure if and only
if it is isomorphic to $\frh_1$ or~$\frh_3$. Moreover, a complex
structure on $\frh_3$ has a compatible LCK metric if and only if
it is equivalent to~$J_0^+$.
\end{theorem}

\begin{proof}
Since the Lee form $\theta$ is a real 1-form, there exist
$\lambda_j\in\C$, $j=1,2,3$, such that
\begin{equation}\label{Leeform}
\theta=\lambda_1\,\omega^1+\lambda_2\,\omega^2+\lambda_3\,\omega^3
+\bar{\lambda}_1\,\omega^{\bar{1}}
+\bar{\lambda}_2\,\omega^{\bar{2}}
+\bar{\lambda}_3\,\omega^{\bar{3}},
\end{equation}
with respect to any basis $\{ \omega^j \}_{j=1}^3$ for
$\frg^{1,0}$. We must find the possible values of $\lambda_j$
in~(\ref{Leeform}) satisfying the equation $\partial\Omega
=\theta^{1,0}\wedge\Omega$. From~(\ref{2forma}) it follows that
\begin{equation}\label{componente21}
\begin{array}{rcl}
\theta^{1,0}\wedge\Omega=&\!\!\! ( \!\!\!&\!\!\!
\lambda_1\,\omega^1+\lambda_2\,\omega^2+\lambda_3\,\omega^3)\wedge\Omega\\[2pt]
=\!\!&\!\!\!-\!\!\!&\!\! (\bar{u}\lambda_1+ir\lambda_2)
\,\omega^{12\bar{1}} +
(is\lambda_1-u\lambda_2)\,\omega^{12\bar{2}} + (v\lambda_1-z\lambda_2)\,\omega^{12\bar{3}}\\
\!\!&\!\!\!-\!\!\!&\!\!
(\bar{z}\lambda_1+ir\lambda_3)\,\omega^{13\bar{1}} -
(\bar{v}\lambda_1+u\lambda_3)\,\omega^{13\bar{2}} +
(it\lambda_1-z\lambda_3)\,\omega^{13\bar{3}}\\
\!\!&\!\!\!-\!\!\!&\!\!
(\bar{z}\lambda_2-\bar{u}\lambda_3)\,\omega^{23\bar{1}} -
(\bar{v}\lambda_2+is\lambda_3)\,\omega^{23\bar{2}} +
(it\lambda_2-v\lambda_3)\,\omega^{23\bar{3}}.
\end{array}
\end{equation}
We shall also use the fact that the closedness of $\theta$ is
equivalent to
$\partial\theta^{1,0}=\db\theta^{1,0}+\partial\theta^{0,1}=0$.

By Lemma~\ref{LCKinvariance} and Proposition~\ref{equivspecial} we
can restrict our attention to the two special types of complex
structures defined by~(\ref{nonnilpotent}) and~(\ref{nilpotent}).
If $J$ is a nonnilpotent complex structure defined
by~(\ref{nonnilpotent}) then $0=\partial\theta^{1,0}=\lambda_2\,
E\, \omega^{13}$, which implies $\lambda_2=0$. Moreover, comparing
the coefficients of $\omega^{23\bar{2}}$ in
Lemma~\ref{partialOmega}~(i) and (\ref{componente21}) we get that
$is \lambda_3 = 0$, so $\lambda_3=0$ because $g$ is positive
definite. Now, if we compare the coefficients of
$\omega^{23\bar{1}}$ then $is-bt=0$, which is a contradiction to
the fact that $s,b,t$ are nonzero real numbers. Therefore, a
nonnilpotent complex structure cannot have compatible LCK metrics.

Let us suppose next that $J$ is a nilpotent complex structure
defined by~(\ref{nilpotent}).

Notice that if the coefficient $\lambda_3$ in~(\ref{Leeform})
vanishes, then comparing the coefficients of $\omega^{13\bar{3}}$
and $\omega^{23\bar{3}}$ in Lemma~\ref{partialOmega}~(ii) and
(\ref{componente21}) we get that $\lambda_1=\lambda_2=0$, so
$d\Omega=0$ and $\frg$ must be the abelian Lie algebra
$\frh_1$~\cite{BG}.

On the other hand, if $\epsilon=1$ in equations~(\ref{nilpotent})
then the coefficients of $\omega^{23\bar{2}}$ and
$\omega^{23\bar{3}}$ in Lemma~\ref{partialOmega}~(ii) and
(\ref{componente21}) imply that $\lambda_2$ and $\lambda_3$
satisfy $\bar{v} \lambda_2 + is \lambda_3 =
 i t \lambda_2 - v \lambda_3 = 0$.
Since $g$ is positive definite, $\det \left(
\begin{array}{lr}
\bar{v} & is \\
it & -v
\end{array} \right)
>0$ and the unique solution is the trivial one,
in particular $\lambda_3=0$ and so $\frg\cong \frh_1$ again.

Suppose next that the NLA $\frg$ is not abelian and it is endowed
with a nilpotent complex structure $J$ given by~(\ref{nilpotent})
admitting an LCK metric. From the previous paragraphs,
$\epsilon=0$ in~(\ref{nilpotent}) and $\lambda_3\not=0$
in~(\ref{Leeform}). From~(\ref{nilpotent}) we have
$0=\partial\theta^{1,0}=\lambda_3\, \rho\, \omega^{12}$, therefore
$\rho=0$ and the complex structure $J$ must be abelian. Since
$\epsilon=\rho=0$, Proposition~\ref{epsilon} and
Lemma~\ref{epsilonzero-new} imply that we can suppose $J$ given by
equations~(\ref{epsilonzero-red-new}) with $\rho=0$. But in this
case one has
$$0=\db\theta^{1,0}+\partial\theta^{0,1}=(\lambda_3-\bar{\lambda}_3)\omega^{1\bar{1}}
+ B\lambda_3\,\omega^{1\bar{2}} -
\bar{B}\bar{\lambda}_3\,\omega^{2\bar{1}} +
(D\lambda_3-\bar{D}\bar{\lambda}_3)\omega^{2\bar{2}}.
$$
In order to have a solution  with $\lambda_3\not=0$, the
coefficient $B$ must be zero and the coefficients $\lambda_3$ and
$D$ must be real. In this case, we get
$$
\partial\Omega = v\, \omega^{12\bar{1}} + D z\, \omega^{12\bar{2}}
-it\, \omega^{13\bar{1}} -iDt\, \omega^{23\bar{2}}.
$$
Now taking into account the coefficients of $\omega^{13\bar{1}}$
and $\omega^{13\bar{3}}$ in~(\ref{componente21}), the condition
$\partial\Omega = \theta^{1,0}\wedge\Omega$ implies that
$\bar{z}\lambda_1+ir\lambda_3=it$ and $it\lambda_1-z\lambda_3 =0$,
so $\lambda_3=t^2/(rt-|z|^2)$. Moreover, from the coefficients of
$\omega^{23\bar{2}}$ and $\omega^{23\bar{3}}$
in~(\ref{componente21}) we get that $\lambda_3 = D\,
t^2/(st-|v|^2)$. Since $g$ is positive definite, necessarily
$D>0$. Now, Corollary~\ref{J0+-} implies that $\frg\cong\frh_3$
and the complex structure $J$ must be equivalent to $J_0^+$.

Finally, the existence of a particular LCK structure on $\frh_3$
follows from~\cite{CFL}. In fact, one solution is obtained for
$D=1$ and $r=s=t=1$, $u=v=z=0$, with Lee form $\theta=2\, \Re
\omega^3$.
\end{proof}

\begin{remark}\label{interseccion}
{\em According to~\cite{BG}, ${\cal M}_3(\frg,J)\cap {\cal
M}_4(\frg,J)=\emptyset$ for any complex structure $J$ on a
nonabelian NLA $\frg$. From Theorems~\ref{balanced-clasif}
and~\ref{LCKgen} we have that for any $J$ on the Lie algebra
$\frh_3$, either ${\cal M}_3(\frh_3,J)=\emptyset$ or ${\cal
M}_4(\frh_3,J)=\emptyset$, depending on the fact that $J$ be
equivalent to $J_0^+$ or not. Moreover, for the remaining
(nonabelian) Lie algebras $\frg$ of Theorem~\ref{clasif-complex},
one has that ${\cal M}_4(\frg,J)=\emptyset$ for any complex
structure $J$.}
\end{remark}

Next we prove that the Lee form of any invariant LCK structure is
parallel with respect to the Levi-Civita connection.

\begin{proposition}\label{genHopfmetrics}
Any invariant LCK metric on the nilmanifold $N(2,1)\times S^1$ is
a generalized Hopf metric.
\end{proposition}

\begin{proof}
Since the complex structure must be equivalent to $J_0^+$, we
consider a basis $\{\omega^j\}_{j=1}^3$ for $(\frh_3)^{1,0}$
satisfying $d\omega^1 = d\omega^2 = 0$ and $d\omega^3 =
\omega^{1\bar{1}}+ \omega^{2\bar{2}}$. It is easy to see that a
$J_0^+$-Hermitian metric $g$ given by~(\ref{metric}) is LCK if and
only if $u=(i\bar{v}z)/ t$ and $|v|^2-st=|z|^2-rt$. In this case,
the associated Lee form is
$$
\theta= {1\over |z|^2-rt} (itz\,\omega^1+itv\,\omega^2-t^2\omega^3
-it\bar{z}\,\omega^{\bar{1}} -it\bar{v}\,\omega^{\bar{2}}
-t^2\omega^{\bar{3}}).
$$
Let $\{Z_j\}_{j=1}^3$ be the dual basis of $\{\omega^j\}_{j=1}^3$.
For any $U,V\in (\frh_3)_\peqC$, it is easy to check that
$$
\theta(\nabla_V U)={it\over |z|^2-rt}\, g(\nabla_V
U,Z_3+\bar{Z}_3).
$$
But, $ g(\nabla_{Z_k} Z_j,Z_3+\bar{Z}_3)=g(\nabla_{Z_k}
\bar{Z}_j,Z_3+\bar{Z}_3)=0$ for $1\leq j\leq k\leq 3$, because
$[Z_1,\bar{Z}_1]=[Z_2,\bar{Z}_2]=\bar{Z}_3-Z_3$ are the only
brackets which do not vanish. Therefore, $g(\nabla_V
U,Z_3+\bar{Z}_3)$ vanishes identically, so the Lee form $\theta$
is parallel.
\end{proof}

It is well-known that the orthogonal leaves to the Lee vector
field of a generalized Hopf manifold bear a Sasakian structure,
and that the product by $\R$ or $S^1$ of a Sasakian manifold is an
LCK manifold with parallel Lee form~\cite{V2}. Thus, as an
immediate consequence of Proposition~\ref{genHopfmetrics} we have
that $N(2,1)$ is essentially the only $5$-dimensional nilmanifold
admitting invariant Sasakian structures:

\begin{corollary}\label{sasakian}
Let $M=\nilm$ be a non-toral compact nilmanifold of dimension $5$
endowed with an invariant Sasakian structure. Then, the Lie
algebra of $G$ is isomorphic to $(0,0,0,0,12+34)$.
\end{corollary}

\medskip

Following an idea of~\cite{Belg}, next we study the symmetrization
of LCK structures on nilmanifolds.

\begin{proposition}\label{symmetrization-LCK}
Let $(M=\nilm,J)$ be a compact complex nilmanifold with $J$
invariant. If $(J,g)$ is an LCK structure on $M$ then there is a
metric $\widetilde{g}$ globally conformal to $g$ such that
$(J,\widetilde{g}_{\nu})$ is an LCK structure on the Lie algebra
$\frg$ of $G$.
\end{proposition}

\begin{proof}
The fundamental form $\Omega$ of $(J,g)$ satisfies
$d\Omega=\theta\wedge\Omega$ with closed Lee form $\theta$. From
Remark~\ref{nomizu} we have that $\theta$ is cohomologous to the
invariant 1-form $\theta_{\nu}$ obtained by the symmetrization
process. Thus, there exists a function $f$ on $M$ such that
$\theta_{\nu}-\theta=df$. Since
$$
d(\exp f\, \Omega)= \exp f\, df\wedge\Omega  + \exp
f(\theta_{\nu}-df)\Omega = \theta_{\nu}\wedge(\exp f\,\Omega),
$$
there is an LCK structure $(J,\, \widetilde{g}=\exp f\, g)$ on $M$
with fundamental form $\widetilde{\Omega}=\exp f\,\Omega$ and Lee
form equal to $\theta_{\nu}$. Now, $d \widetilde{\Omega}_{\nu} =
(d\widetilde{\Omega})_{\nu} = (\theta_{\nu}\wedge
\widetilde{\Omega})_{\nu} = \theta_{\nu}\wedge
\widetilde{\Omega}_{\nu}$, that is, $(J,\widetilde{g}_\nu)$ is an
LCK structure on the Lie algebra~$\frg$.
\end{proof}

\begin{corollary}\label{LCK-general}
Let $(M=\nilm,J)$ be a non-toral $6$-dimensional compact
nilmanifold endowed with an invariant complex structure $J$. Then,
$M$ has an LCK metric if and only if the Lie algebra of $G$ is
isomorphic to $\frh_3$ and $J$ is equivalent to $J_0^+$.
\end{corollary}

It is a conjecture of Vaisman that any compact locally conformal
K\"ahler but not globally conformal K\"ahler manifold has an odd
Betti number. By Corollary~\ref{LCK-general} this conjecture is
true in the class of compact nilmanifolds with invariant complex
structure up to dimension six. In this context, it seems natural
to conjecture that {\it a $2n$-dimensional compact nilmanifold $M$
admitting LCK structure is the product of $N(n-1,1)$ by $S^1$},
where $N(n-1,1)$ is the quotient of the generalized Heisenberg
group $H(n-1,1)$ by a discrete subgroup, in particular the first
Betti number of $M$ equals $2n-1$; that is to say, the only LCK
nilmanifolds are essentially those constructed in~\cite{CFL}.

\medskip

The following result shows a large class of complex nilmanifolds
not admitting LCK structures.

\begin{corollary}\label{LCK-complex-par}
A compact complex parallelizable nilmanifold (not a torus) has no
LCK metrics.
\end{corollary}

\begin{proof}
Let $M$ be a compact complex parallelizable nilmanifold and denote
by $J$ its complex structure. Since $M$ is not a torus and any
invariant $J$-Hermitian metric is balanced~\cite{AG}, there do not
exist invariant LCK metrics on $M$. By
Proposition~\ref{symmetrization-LCK} there are no LCK metrics on
$M$ compatible with $J$.
\end{proof}

In~\cite[Remark 1]{AI} it is proved that the SKT condition is
complementary to the LCK condition. Next we give another proof of
this fact for nilmanifolds, based on the nilpotency of the
underlying Lie algebra.

\begin{proposition}\label{SKT-LCK}
Let $(M=\nilm,J)$ be a non-toral compact complex nilmanifold of
dimension $2n \geq 6$, where $J$ is invariant. A $J$-Hermitian
metric $g$ on $M$ cannot be SKT and LCK at the same time.
\end{proposition}

\begin{proof}
Let $(J,g)$ be a Hermitian structure on $M$ that is both SKT and
LCK. From Propositions~\ref{symmetrization-SKT}
and~\ref{symmetrization-LCK}, there is a Hermitian structure on
the Lie algebra $\frg$ of $G$ that is SKT and LCK at the same
time, i.e. its fundamental form $\Omega$ satisfies
$d\Omega=\theta\wedge\Omega$ and $\partial\db\Omega=0$. Let us
write the Lee form as $\theta=\theta^{1,0}+\theta^{0,1}$, where
$\theta^{0,1}=\overline{\theta^{1,0}}$. Since
$\theta^{1,0}\wedge\Omega=\partial\Omega$ and
$\db(\theta^{1,0}\wedge\Omega)=-\partial\db\Omega=0$, we have that
$\theta^{1,0}\wedge\Omega$ is a closed form. Therefore,
$0=d(\theta^{1,0}\wedge\Omega)=(d\theta^{1,0}-\theta^{1,0}\wedge\theta^{0,1})\wedge\Omega$,
which implies $d\theta^{1,0}=\theta^{1,0}\wedge\theta^{0,1}$,
because the dimension of $\frg$ is $\geq 6$. Notice that
$\theta^{1,0}\not= 0$ because $\frg$ is not abelian. Now, the real
1-form $\eta=i(\theta^{1,0}-\theta^{0,1})$ satisfies
$d\eta=\eta\wedge\theta$, and a standard argument shows that this
cannot happen because $\frg$ is nilpotent.
\end{proof}

\begin{remark}\label{dim-4}
{\rm The proposition above does not hold for nilmanifolds of
dimension 4. In fact, for any complex structure on the Lie algebra
$\frk\frt=(0,0,0,12)$ underlying the well-known Kodaira-Thurston
manifold~\cite{Th}, there is a basis $\{\omega^1,\omega^2\}$ of
$\frk\frt^{1,0}$ such that $d\omega^1=0$ and
$d\omega^2=\omega^{1\bar{1}}$. For any compatible metric
$$
g= r\, \omega^1\#\omega^{\bar{1}} + s\,\omega^2\#\omega^{\bar{2}}
 - i( u\,\omega^1\#\omega^{\bar{2}}
-\bar{u}\,\omega^2\#\omega^{\bar{1}} ),
$$
its fundamental form $\Omega$ satisfies $\partial\db\Omega=0$, so
$g$ is SKT. Moreover, $g$ is also LCK, because
$d\Omega=\theta\wedge\Omega$ with closed
$\theta=\frac{2s}{|u|^2-rs}(\Re(iu\omega^1) - s\, \Re \omega^2)$.}
\end{remark}

\bigskip

\noindent {\bf Acknowledgments}.- This work has been partially
supported by project BFM2001-3778-C03-03.

{\small

\vspace{0.25cm}

\noindent{\sf Departamento de Matem\'aticas, Facultad de Ciencias,
Universidad de Zaragoza, Campus Plaza San Francisco, 50009
Zaragoza, Spain.} {\sl E-mail:} ugarte@unizar.es }


\begin{thebibliography}{33}

\bibitem{AG} E. Abbena, A. Grassi, Hermitian left invariant
  metrics on complex Lie groups and cosymplectic Hermitian
  manifolds, {\it Boll. U.M.I.\/} {\bf A} (6) {\bf 5} (1986),
  371--379.

\bibitem{AI} B. Alexandrov, S. Ivanov, Vanishing theorems on Hermitian manifolds,
 {\it Differ. Geom. Appl.\/} {\bf 14} (2001), 251--265.

\bibitem{Belg}
  F.A. Belgun, On the metric structure of non-K\"ahler complex surfaces,
  {\it Math. Ann.\/} {\bf 317} (2000), 1--40.

\bibitem{BG}
  C. Benson, C.S. Gordon, K\"ahler and symplectic structures on
  nilmanifolds, {\it Topology\/} {\bf 27} (1988), 513--518.

\bibitem{Bismut}
  J.-M. Bismut, A local index theorem for non-K\"ahler manifolds,
  {\it Math. Ann.\/} {\bf 284} (1989), 681--699.

\bibitem{palermo}
  L.A. Cordero, M. Fern\'andez, A. Gray, L. Ugarte, Nilpotent
  complex structures on compact nilmanifolds, {\it Rend. Circolo
  Mat. Palermo} {\bf 49} {\it suppl.} (1997), 83--100.

\bibitem{CFGU2}
  L.A. Cordero, M. Fern\'andez, A. Gray, L. Ugarte,
  Compact nilmanifolds with nilpotent complex structure:
  Dolbeault cohomology,
  {\it Trans. Amer. Math. Soc.\/} {\bf 352} (2000), 5405--5433.

\bibitem{CFL}
  L.A. Cordero, M. Fern\'andez, M. de Le\'on,
  Compact locally conformal K\"{a}hler nilmanifolds,
  {\it Geom. Dedicata\/} {\bf 21} (1986), 187--192.

\bibitem{abelian}
  L.A. Cordero, M. Fern\'andez, L. Ugarte,
  Abelian complex structures on 6-dimensional compact nilmanifolds,
  {\it Comment. Math. Univ. Carolinae\/} {\bf 43} (2002), 215--229.

\bibitem{pseudoKahler}
  L.A. Cordero, M. Fern\'andez, L. Ugarte,
  Pseudo-K\"ahler metrics on six-dimensional nilpotent Lie algebras,
  {\it J. Geom. Phys.\/} {\bf 50} (2004), 115--137.

\bibitem{DO}
  S. Dragomir, L. Ornea, {\sl Locally Conformal K\"ahler
  Geometry}, Progress in Math. 155, Birkh\"auser, 1998.

\bibitem{FG}
  A. Fino, G. Grantcharov, Properties of manifolds with
  skew-symmetric torsion and special holonomy, {\it Adv. Math.\/}
  {\bf 189} (2004), 439--450.

\bibitem{FPS}
  A. Fino, M. Parton, S. Salamon, Families of strong KT structures
  in six dimensions, {\it Comment. Math. Helv.\/}
  {\bf 79} (2004), 317--340.

\bibitem{Gau}
  P. Gauduchon, Hermitian connections and Dirac operators, {\it Boll. Un. Mat. Ital.\/}
  {\bf B 11} (1997), 257--288.

\bibitem{GH}
  A. Gray, L.M. Hervella,
  The sixteen classes of almost Hermitian manifolds and their linear
  invariants,
  {\it Ann. Mat. Pura Appl.\/} (4) {\bf 123} (1980), 35--58.

\bibitem{GIP}
  J. Gutowski, S. Ivanov, G. Papadopoulos,
  Deformations of generalized calibrations and compact non-K\"ahler
  manifolds with vanishing first Chern class,
  {\it Asian J. Math.\/} {\bf 7} (2003), 39--80.

\bibitem{L}
  J. Lauret, Geometric structures on nilpotent Lie groups: on their classification
  and a distinguished compatible metric, preprint DG/0210143.

\bibitem{Mal}
  I.A. Mal'cev, A class of homogeneous spaces, {\it Amer. Math.
  Soc. Transl.\/} No. {\bf 39} (1951).

\bibitem{MPPS}
  C. McLaughlin, H. Pedersen, Y.S. Poon, S. Salamon, Deformation of 2-step nilmanifolds
  with abelian complex structures, preprint DG/0402069.

\bibitem{Milnor}
  J. Milnor, Curvature of left invariant metrics on Lie groups,
  {\it Adv. Math.\/} {\bf 21} (1976), 293--329.

\bibitem{No}
  K. Nomizu, On the cohomology of compact homogeneous spaces of
  nilpotent Lie groups, {\it Ann. of Math.\/} {\bf 59} (1954), 531--538.

\bibitem{P}
  G. Papadopoulos,
  KT and HKT geometries in strings and in black hole moduli spaces,
  preprint hep-th/0201111.

\bibitem{S}
  S. Salamon, Complex structures on nilpotent Lie algebras,
 {\it J. Pure Appl. Algebra\/} {\bf 157} (2001), 311--333.

\bibitem{Th}
  W.P. Thurston, Some simple examples of symplectic
  manifolds, {\it Proc. Amer. Math. Soc.\/} {\bf 55} (1976), 467--468.

\bibitem{V2}
  I. Vaisman,
  Locally conformal K\"ahler manifolds with parallel Lee form,
  {\it Rend. Mat.\/} (2) {\bf 12} (1979), 263--284.

\bibitem{V3}
  I. Vaisman,
  Generalized Hopf manifolds,
  {\it Geom. Dedicata\/} {\bf 13} (1982), 231--255.

\bibitem{Wang}
  H.C. Wang,
  Complex parallelizable manifolds,
  {\it Proc. Amer. Math. Soc.\/} {\bf 5} (1954), 771--776.

\end{thebibliography}
\end{document}